\documentclass[12pt]{amsart}
\usepackage{graphicx}
\usepackage{verbatim}
\usepackage{amsmath}
\usepackage{amssymb, mathrsfs}
\usepackage{amsbsy}
\usepackage{amscd}
\usepackage{amsthm}
\usepackage{bm}
\newcommand{\C}{\mathbb{C}}

\newcommand{\N}{\mathbb{N}}
\newcommand{\R}{\mathbb{R}}
\newcommand{\Z}{\mathbb{Z}}
\newcommand{\E}{\mathbb{E}}
\newcommand{\set}[1]{\left\{#1\right\}}
\newcommand{\e}{\epsilon}

\renewcommand{\P}{\mathbb{P}}

\newcommand{\I}{\mathcal{I}}
\renewcommand{\d}{{\rm{d}}}

\newtheorem{theorem}{Theorem}
\newtheorem*{theorem*}{Theorem}
\newtheorem{lemma}{Lemma}
\newtheorem{proposition}{Proposition}
\newtheorem*{proposition*}{Proposition}

\theoremstyle{definition}
\newtheorem{remark}{Remark}
\newtheorem{remarks}[remark]{Remarks}

\numberwithin{equation}{section}
\begin{document}
\title[Brownian motion in a heavy tailed Poissonian potential]
{Second order asymptotics for Brownian motion in a heavy tailed Poissonian potential}
\author{Ryoki Fukushima}
\thanks{
This work was partially supported by JSPS Fellowships for Young Scientists and 
a grant of the Swiss National Foundation No.~200020-125247/1.}
\address{
{\rm Ryoki Fukushima}\\
Department of Mathematics, Kyoto University, Kyoto 606-8502, JAPAN}
\curraddr{Department of Mathematics, Tokyo Institute of Technology, Tokyo 152-8551, JAPAN}
\email{ryoki@math.titech.ac.jp}
\keywords{Brownian motion; random media; parabolic Anderson model; Poissonian potential}
\subjclass[2010]{Primary 60K37; Secondary 82B44}

\begin{abstract}
We consider the Feynman-Kac functional associated with a Brownian motion in a random potential.
The potential is defined by attaching 
a heavy tailed positive potential around a Poisson point process. 
This model was first considered by Pastur~[Teoret. Mat. Fiz. 32(1), 88--95 (1977)] 
and the first order term of the moment asymptotics was determined. 
In this paper, both moment and almost sure asymptotics are determined up to the second order. 
We also derive the second order asymptotics of the integrated density of states 
of the corresponding random Schr\"{o}dinger operator. \\
\keywords{Brownian motion \and parabolic Anderson model \and random media \and Poissonian potential}
\end{abstract}

\maketitle

\section{Introduction}\label{Intro}
We consider a Brownian motion moving in a Poissonian potential. 
Such a process can be viewed as a polymer in a random environment and 
also has links to other topics in random media such as 
the spectral theory of random Schr\"{o}dinger operator or 
intermittency for a parabolic problem with random potential. 
We refer to a review article~\cite{GK04} by G\"artner and K\"onig
for background and Sznitman's monograph~\cite{Szn98} for a thorough 
study on a model similar to ours. 

Let $(\{B_t\}_{t \ge 0}, P_0)$ be a Brownian motion on $\R^d$ 
with generator $-\kappa \Delta$, 
starting at the origin. 
We define the random potential by attaching the shape function 
$\hat{v}(x)=|x|^{-\alpha}\wedge 1$
around a Poisson point process $(\omega=\sum_i \delta_{\omega_i}, \P_{\nu})$ with constant intensity $\nu$ 
as follows: 
\begin{equation*}
	V_{\omega}(x)=\sum_i \hat{v}(x-\omega_i). 
\end{equation*}
We shall also use notation $v(x)=|x|^{-\alpha}$ later. 
The main object in this work is the long time asymptotics of the Feynman-Kac functional
\begin{equation}
	E_0 \left[ \exp\set{-\int_0^t V_{\omega}(B_s) {\rm d} s} \right].\label{F-K}
\end{equation}
It is well known that \eqref{F-K} represents the total mass of the minimal solution 
of the initial value problem
\begin{equation}
\begin{aligned}
	\frac{\partial}{\partial t}v(t,x) &= \kappa\Delta v(t,x) - V_{\omega}(x) v(t,x), \quad 
	&&(t, x) \in (0, \infty) \times \R^d, \\
	v(0,x) &= \delta_0(x), && x \in \R^d. \label{PAM}
\end{aligned}
\end{equation}

One can also identify \eqref{F-K} as the survival probability of the
Brownian motion killed by the 
random potential $V_{\omega}$. We refer the reader to Section 1.3 of~\cite{Szn98} 
for an illustrative construction of such a process.

\subsection{Early studies}
We mention some early studies which are related to ours. 
When $\alpha > d+2$, which is referred to as the light tailed case, 
Donsker and Varadhan~\cite{DV75c} determined the moment asymptotics 
\begin{equation}
\begin{split}
	&\E_{\nu} \otimes E_0 \left[ \exp\set{-\int_0^t V_{\omega}(B_s) {\rm d} s} \right]\\
	&\quad= \exp\set{-\frac{d+2}{2}(\nu \omega_d)^{\frac{2}{d+2}}
	 \left( \frac{\kappa \lambda_d}{d} \right)^{\frac{d}{d+2}}t^{\frac{d}{d+2}}
	 (1+o(1))}\label{DV}
\end{split}
\end{equation}
as $t$ goes to $\infty$, where $\omega_d$ is the volume of the unit ball and $\lambda_d$ 
the smallest Dirichlet eigenvalue 
of $-\Delta$ in such a ball. Note that the leading asymptotics depends on $\kappa$ but 
not on $\alpha$. In fact, the constant in front of $t^{d/(d+2)}$ has a
variational 
expression
\begin{equation}
	\inf_{U: \,{\rm open}}\set{\kappa \lambda_1(U) + \nu |U|}\label{DV-var}
\end{equation}
for any $\alpha > d+2$, where $|U|$ and $\lambda_1(U)$ stand for the
volume of 
$U$ and the smallest Dirichlet eigenvalue of 
$-\Delta$ in $U$, respectively.
It follows from Faber-Krahn's inequality that balls with radius
\begin{equation}
 R_0=\left(\frac{\kappa \lambda_d}{d\nu\omega_d}\right)
\end{equation}
are the minimizers of~\eqref{DV-var}. 
Later, Sznitman~\cite{Szn93c} proved that when 
$\hat{v}$ has compact support, 
\begin{equation}
 E_0 \left[ \exp\set{-\int_0^t V_{\omega}(B_s) {\rm d} s} \right]
 = \exp\set{-\kappa\lambda_d 
 \left(\frac{\nu \omega_d}{d}\right)^{\frac{2}{d}}
 t (\log t)^{-\frac{2}{d}}
 (1+o(1))}\label{LT-q}
\end{equation}
as $t$ goes to $\infty$, $\P_{\nu}$-almost surely. This case can be
thought of as $\alpha=\infty$. 

On the other hand, in the heavy tailed case $d < \alpha < d+2$, Pastur~\cite{Pas77} 
determined the moment asymptotics
\begin{equation}
	\E_{\nu} \otimes E_0 \left[ \exp\set{-\int_0^t V_{\omega}(B_s) {\rm d} s} \right]
	= \exp\set{-(a_1+o(1))t^{\frac{d}{\alpha}}}\label{Pas}
\end{equation}
as $t$ goes to $\infty$, where 
\begin{equation*}
	a_1=\nu \omega_d \Gamma\left( \frac{\alpha-d}{\alpha} \right). 
\end{equation*}
In contrast to \eqref{DV}, this leading term does not depend on 
$\kappa$ and is determined only by the potential. Indeed, a key step in Pastur's proof 
was proving the asymptotic equivalence between the left hand side of \eqref{Pas} and 
$\E_{\nu}[\exp\{ -t V_{\omega}(0) \}]$.

\begin{remarks}
\begin{enumerate}
 \item For the critical case $\alpha=d+2$, 
 we refer the interested reader to \^{O}kura~\cite{Oku81}. 
 \item Chen and Kulik~\cite{CK11a} has recently proved that for some 
class of potentials $K$, including $|x|^{-\alpha}$ $(\alpha>d/2)$, 
the ``renormalized Poisson potential'' 
\begin{equation*}
 \bar{V}_{\omega}(x)=\int K(x-y)(\omega(\d x)-\nu\d x)
\end{equation*}
can be properly defined and associated parabolic Anderson problem
admits a Feynman-Kac solution. 
In the subsequent papers by Chen {\em{et al.}}~\cite{C11,CK11b,CR11}, 
both moment and almost sure asymptotics of the Feynman-Kac sulution
have been investigated for the case $K(x)=|x|^{-\alpha}$ ($d/2<\alpha<d$).
\end{enumerate}
\end{remarks}
\subsection{Motivation and results}
We shall mainly discuss the second order asymptotics of \eqref{Pas} and
the almost sure asymptotics up to the second term in the heavy tailed case. 
Let us briefly explain why we are interested in higher order terms.

In the light tailed case, Donsker and Varadhan's 
result suggests that the dominant contribution to the right hand side 
of~\eqref{DV} comes from the following strategy: 
there exists $x\in\R^d$ such that $\omega(B(x,R_0t^{1/(d+2)}))=0$ and
the Brownian motion $\{B_s\}_{0\le s\le t}$ stays in the ball. 
Indeed, one can easily see that this specific event gives the 
correct lower bound. 
Motivated by this observation, 
Sznitman~\cite{Szn91b} ($d=2$) and Povel~\cite{Pov99} ($d \ge 3$) 
proved that the above confinement is typical for the paths which
survives for a long time. 
Sznitman also studied the behavior of surviving paths under a 
fixed configuration $\omega$, based on the heuristics behind
\eqref{LT-q}.  
We refer the reader to the monograph~\cite{Szn98} for detail and 
related topics. 

Our motivation is to develop the study on the typical behavior of the 
surviving paths in the heavy tailed case. 
However in the heavy tailed case, it seems difficult to read
so much 
information about 
the Brownian motion from \eqref{Pas} since it is independent of the
diffusion 
coefficient $\kappa$. 

The first main theorem of this article is the moment asymptotics of the 
Feynman-Kac functional 
up to the second order. 
It in particular gives refinement of \eqref{Pas} and we see that the
second 
term does depend on $\kappa$. 
\begin{theorem}\label{annealed}
Suppose $d < \alpha < d+2$ and $p \in [0, \infty)$. Then
\begin{equation}
\begin{split}
	&\E_{\nu} \left[ E_0 \left[\exp\set{-\int_0^t 
        V_{\omega}(B_s) {\rm d} s} \right]^p \right]\\
	&\quad = \exp\set{-a_1(pt)^{\frac{d}{\alpha}} -
        (a_2+o(1))(pt)^{\frac{\alpha+d-2}{2\alpha}}}\label{thm1}
\end{split}
\end{equation}
as $t \to \infty$, where
\begin{equation*}
	a_2 
	= \left(\frac{\kappa \nu \alpha \sigma_d}{2} 
        \Gamma\left(\frac{2\alpha-d+2}{\alpha}\right)\right)^{\frac{1}{2}}.
\end{equation*}
Moreover, the constant $a_2$ admits a variational expression
\begin{equation}
 a_2=\inf_{\stackrel{\phi \in W^{1,2}(\R^d), }{ \| \phi \|_{L^2}=1}}
	\set{\int \kappa|\nabla \phi |(x)^2 +
        C(\nu,d,\alpha)
	|x|^2\phi(x)^2{\rm d} x }\label{var-Fuk}\\
\end{equation}
with
\begin{equation*}
 C(\nu, d, \alpha)=\frac{\nu\alpha\sigma_d}{2d} 
        \Gamma\left( \frac{2\alpha-d+2}{\alpha} \right),
\end{equation*}
where $\sigma_d$ denotes the surface area of the unit sphere and 
$W^{1,2}(\R^d)$ the usual Sobolev space. 
\end{theorem}

The second main result is the almost sure asymptotics of the Feynman-Kac functional. 
The dependence on the diffusion coefficient appears in the second term again. 
\begin{theorem}\label{quenched}
Suppose $d < \alpha < d+2$. Then for $\P_{\nu}$-almost every $\omega$, 
\begin{equation*}
\begin{split}
	&E_0 \left[ \exp\set{-\int_0^t V_{\omega}(B_s) {\rm d} s} \right]\\
	&\quad= \exp\set{-q_1 t(\log t)^{-\frac{\alpha-d}{d}}
	-\left(q_2 +o(1)\right)t(\log t)^{-\frac{\alpha-d+2}{2d}}}
\end{split}
\end{equation*}
as $t \to \infty$, where
\begin{equation*}
\begin{split}
	& q_1= \frac{d}{\alpha} \left(\frac{\alpha-d}{\alpha d} 
        \right)^{\frac{\alpha-d}{d}}
	a_1^{\frac{\alpha}{d}},\\
	& q_2= a_2\left(\frac{\alpha-d}{\alpha d} a_1\right)^{\frac{\alpha-d+2}{2d}}.\\
\end{split}
\end{equation*}
\end{theorem}

Finally we state our result on 
the integrated density of states of 
the random Schr\"{o}dinger operator $-\kappa\Delta+V_{\omega}$. 
Recall that the integrated density of states is defined by  
\begin{equation}
	N(\lambda) = \lim_{R \to \infty} \frac{1}{(2R)^d}\E_{\nu}
	\bigl[\#\bigl\{ k \in \N ;\lambda_{\omega,\,k}\bigl((-R, R)^d\bigr) \le \lambda \bigr\}\bigr],
	\label{IDS-def} 
\end{equation}
where $\lambda_{\omega,\,k}\bigl((-R, R)^d\bigr)$ is the 
$k$-th smallest eigenvalue of $-\kappa\Delta+V_{\omega}$ in $\bigl((-R, R)^d\bigr)$
with the Dirichlet boundary condition. The existence of the above limit is proved for instance in~\cite{KM82a}. 
\begin{remark}\label{Rem-IDS}
The following statements are also proved in~\cite{KM82a}: 
\begin{enumerate}
\item{The limit without taking $\E_{\nu}$ in \eqref{IDS-def} exists almost surely
and coincides with \eqref{IDS-def}. 
This is more usual definition of the integrated density of states.}
\item{In our setting of $V_{\omega} \ge 0$, the limit in \eqref{IDS-def} is unchanged 
if we consider the Neumann boundary condition instead. }
\item{The existence of the limit in \eqref{IDS-def} is proved by using a
     spatial superadditivity property of $\E_{\nu}
	[\#\{ k \in \N ;\lambda_{\omega,\,k}((-R,
     R)^d) 
     \le \lambda \}]$. 
Hence it is in fact the supremum over $R>0$. }
\end{enumerate}
\end{remark}

On the way of the proof of Theorem~\ref{quenched}, we 
obtain the second order asymptotics of the integrated density of states. 
\begin{theorem}\label{Lifshitz}
Suppose $d < \alpha < d+2$. Then
\begin{equation}
	N(\lambda)=\exp\set{ -l_1 \lambda^{-\frac{d}{\alpha-d}} 
	- (l_2+o(1)) \lambda^{-\frac{\alpha+d-2}{2(\alpha-d)}} }\label{Lifshiz}
\end{equation}
as $\lambda \downarrow 0$, where 
\begin{align*}
	&l_1=\frac{\alpha-d}{\alpha} \left(\frac{d}{\alpha} \right)^{\frac{d}{\alpha-d}}
	a_1^{\frac{\alpha}{\alpha-d}},\\
	&l_2=a_2 \left( \frac{da_1}{\alpha} \right)^{\frac{\alpha+d-2}{2(\alpha-d)}}.
\end{align*}
\end{theorem}
The first term in \eqref{Lifshiz} has been determined in~\cite{Pas77} by
using Tauberian theorem and~\eqref{Pas}.
Due to its independent of $\kappa$, it is said that 
the first term of $N(\lambda)$ has classical character 
(see e.g.~\cite{LMW03}).
Our result shows that the quantum effect appears in the second term.

\subsection{Ideas and heuristics}
To understand the ideas and heuristics behind Theorem~\ref{annealed} 
and~\ref{quenched}, it is illustrative to 
see which kind of strategy gives the lower bound, as in the 
light tailed case. 

Let us start with the picture behind 
Pastur's first order asymptotics~\eqref{Pas}. 
The lower bound in~\eqref{Pas} comes from the following strategy: 
\begin{enumerate}
 \item $V_{\omega}(0) \sim a_1\frac{d}{\alpha}t^{-(\alpha-d)/\alpha}$ 
       and
 \item $\{B_s\}_{0\le s \le t}$ is confined in the ball of radius
       $o(t^{1/\alpha})$ centered at 0.
\end{enumerate}
It is not difficult to see that conditioned on (i), 
$V_{\omega}(x) \sim a_1\frac{d}{\alpha}t^{-(\alpha-d)/\alpha}$
for all $|x|=o(t^{1/\alpha})$ with high probability. 
Roughly speaking, this is because $V_{\omega}$ is locally 
stiff where it takes small value. 
Now we can see this strategy indeed gives the correct lower bound
since
\begin{equation}
\begin{split}
  \exp\set{-\int_0^tV_{\omega}(B_s)\d s}
 &\approx \exp\set{-\int_0^tV_{\omega}(0)\d s}\\
 &\approx \exp\set{-a_1\frac{d}{\alpha}t^{\frac d\alpha}}\label{F-K1}
\end{split}
\end{equation}
on this event, the first event (i) has probability 
\begin{equation}
 \P_{\nu}\left(V_{\omega}(0) \sim a_1\frac{d}{\alpha}
 t^{-\frac{\alpha-d}\alpha}\right)
 =\exp\set{-a_1\frac{\alpha-d}{\alpha}t^{\frac d\alpha}(1+o(1))},
\label{probab} 
\end{equation}
and the probability of the second event (ii) is easily 
seen to be $\exp\{-o(t^{d/\alpha})\}$ if the ball is 
not too small (recall that $\alpha<d+2$).

The lower bound of Theorem~\ref{annealed} can be obtained
by a finer analysis of the above strategy. We assume 
for simplicity that 
``$\sim$'' and ``$o(1)$'' in~\eqref{probab} are sufficiently 
precise not to affect the second term in~\eqref{thm1}; 
see Remark~\ref{precision} below. 
The key fact is that conditioned on 
$V_{\omega}(0) \sim a_1\frac{d}{\alpha}t^{-(\alpha-d)/\alpha}$, 
the potential viewed from the bottom locally 
looks like a parabola:
\begin{equation*}
 V_{\omega}(x)-V_{\omega}(0)
 \sim C(\nu,d,\alpha)t^{-\frac{\alpha-d+2}{\alpha}} |x|^2,
\end{equation*}
for $|x|=o(t^{1/\alpha})$. 
Thus we have 
\begin{equation*}
\begin{split}
 &E_0 \left[ \exp\set{-\int_0^tV_{\omega}(B_s)\d s}:
 \sup_{0 \le s \le t}|B_s|=o(t^{\frac{1}{\alpha}})\right]\\
 &\quad \approx \exp\set{-a_1\frac{d}{\alpha}t^{\frac d\alpha}}
 E_0 \biggl[ \exp\set{-Ct^{-\frac{\alpha-d+2}{\alpha}} 
 \int_0^t |B_s|^2 \d s}
 :\sup_{0 \le s \le t}|B_s|=o(t^{\frac{1}{\alpha}})
 \biggr]
\end{split}
\end{equation*}
instead of~\eqref{F-K1}. 
The second term in the right hand side, 
together with the suitable scaling
and Donsker-Varadhan's large deviation theory, explains how 
the term $a_2t^{(\alpha+d-2)/2\alpha}$ arises.
Moreover, it suggests that the typical surviving paths live 
in the scale $t^{(\alpha-d+2)/4\alpha}$ and its (scaled) 
occupation time measure looks like a Gaussian density, 
which is the unique minimizer of~\eqref{var-Fuk}.
Rigorous proof of this heuristics is an interesting 
problem and will be addressed in future work. 

\begin{remark}\label{precision}
It is possible to prove a refinement of~\eqref{probab}
and make the above argument rigorous when $\alpha \ge 2$.
When $1<\alpha<2$, we can follow essentially the same line 
but with some modification.  
We do not go into further detail in this paper since 
the actual proof of Theorem~\ref{annealed} is given by 
applying a abstract theory developed by G\"artner 
and K\"onig~\cite{GK00}. 
\end{remark}

Next we explain the heuristics behind the lower bound of 
Theorem~\ref{quenched}. 
It is natural to expect that the main contribution comes from
paths which spend most of the time in {\em valleys} where 
$V_{\omega}$ takes atypically small value, 
as in~\cite{Szn93c}.
So we first fix a large enough box $(-t,t)^d$ and 
consider the minimum value of $V_{\omega}$ in it. 
Then it follows essentially from~\eqref{probab} that 
the minimum is asymptotic to $q_1(\log t)^{-(\alpha-d)/d}$. 
Furthermore, one can show that if 
$V_{\omega}(m)$ is close to the minimum, 
then the potential locally looks like a parabola around $m$:
\begin{equation}
	q_1(\log t)^{-\frac{\alpha-d}{d}}+
	\frac{q_2^2}{\kappa d}(\log t)^{-\frac{\alpha-d+2}{d}}
        |x-m|^2 \label{valley}
\end{equation}
for $x \in B(m, M(\log t)^{(\alpha-d+2)/2d})$ ($M>0$: large). 
We can obtain the lower bound of Theorem~\ref{quenched} 
by considering the paths which go into one of the valleys 
in relatively short time 
and stay there afterward. 
It is another future problem to show that this strategy is typical for 
surviving paths, which is indeed proved for compactly supported 
$\hat{v}$ by Sznitman~\cite{Szn96B}. 
\begin{remark}
Strictly speaking, \eqref{valley} is proved only for $\alpha \ge 2$
(see Proposition~\ref{alpha>2}).
When $1<\alpha<2$, only a modified version
of~\eqref{valley}, in terms of 
eigenvalue, is proved due to a technical difficulty. 
\end{remark}

\subsection{Outline}
The reminder of the paper is organized as follows.
Theorem~\ref{annealed} is proved in Section~\ref{proof-annealed}. 
The upper bound of 
Theorem~\ref{Lifshitz} is proved in Section~~\ref{Lifshiz-upper}
by applying Tauberian argument to Theorem~\ref{annealed} 
and then used to prove the upper bound of Theorem~\ref{quenched}
at the beginning of Section~\ref{proof-quenched}. 
The lower bound of Theorem~\ref{quenched} is 
the most involved part and constitutes a large portion of
Section~~\ref{proof-quenched}.
It has two subsections since 
we discuss the case $\alpha \ge 2$ and $1<\alpha<2$ separately.
The proof of the lower bound of Theorem~\ref{Lifshitz} 
is given in Section~~\ref{Lifshiz-lower}, using 
an eigenvalue estimate derived in Section~\ref{proof-quenched}. 
Finally, in Section~\ref{appendix} we collect some 
formulae concerning Poisson point process. 

\section{Moment asymptotics}\label{proof-annealed}
For the moment asymptotics, there is a general framework developed by G\"artner 
and K\"onig~\cite{GK00} and we shall make use of it. We first recall the element of their result. 

Let $(\{\xi(x)\}_{x \in \R^d}, \P)$ be a translation invariant random field having 
all positive exponential moments:
\begin{equation*}
	H(t)=\log \E [e^{t \xi(0)}]<\infty \quad \textrm{for}\quad t>0. 
\end{equation*}
For a compact set $K \subset \R^d$, let $\mathcal{P}(K)$ denote the set of probability measures 
whose supports are contained in $K$ and $\mathcal{P}_c(\R^d)$ the set of all compactly supported 
probability measures. 
The main assumption (called Assumption (J)) in~\cite{GK00} is that there exists a 
\emph{scale} $r=r(t)$ such that the functional
\begin{equation*}
	J_t(\mu)=-\frac{1}{tr^{-2}} \left( \log \E\left[ \exp\set{t \int \xi(rx) \,\mu({\rm d} x)}\right] - H(t)\right)
\end{equation*}
on $\mathcal{P}_c(\R^d)$ converges to a functional $J: \mathcal{P}_c(\R^d) \to [0,\infty)$ as $t \to \infty$ 
uniformly on $\mathcal{P}(K)$ for each compact set $K \subset \R^d$. 
For a part of their result, they also require the following (called Assumption (H)): 
\begin{equation*}
	H\left(t+e^{-\e tr^{-2}}\right) - H(t) 
	= O\left(e^{\e tr^{-2}}\right) \quad {\rm as} \quad t \to \infty 
\end{equation*}
for each $\e>0$. 
We need some more notations to state the result. For $\mu \in \mathcal{P}_c(\R^d)$, let 
\begin{equation*}
	I(\mu) =
	\begin{cases}
		\Bigl \| \nabla \sqrt{\frac{{\rm d}\mu}{{\rm d} x}}\Bigr\|_2^2, 
		&\textrm{if } {\rm d}\mu \ll {\rm d} x \textrm{ and } 
                \sqrt{\frac{{\rm d}\mu}{{\rm d} x}}\in W^{1,2}(\R^d), \\
		\infty, & {\rm otherwise} 
	\end{cases}
\end{equation*}
and define 
\begin{equation*}
	\chi = \inf \set{\kappa I(\mu) + J(\mu): \mu \in \mathcal{P}_c(\R^d)}.
\end{equation*}
Then, their main result (Theorem 1 in~\cite{GK00}) is the following. 
\begin{theorem}\label{GK}
Fix $p \in [0, \infty)$ arbitrarily and suppose that {\upshape Assumption (J)} is satisfied. 
\begin{itemize}
\item[\upshape{(i)}]{As $t \to \infty$, 
\begin{equation*}
	\E\left[ E_0 \left[ \exp\set{\int_0^t \xi (B_s) {\rm d} s} \right]^p\right]
	\ge \exp\set{H(pt)-\frac{pt}{r(pt)^2}(\chi+o(1))}.
\end{equation*}
}
\item[\upshape{(ii)}]{If $p=1$ or, in addition,{ \upshape Assumption (H)} 
is satisfied, then, as $t \to \infty$, 
\begin{equation*}
	\E\left[ E_0 \left[ \exp\set{\int_0^t \xi (B_s) {\rm d} s} \right]^p\right]
	\le \exp\set{H(pt)-\frac{pt}{r(pt)^2}(\chi+o(1))}.
\end{equation*}
}
\end{itemize} 
\end{theorem}

From now on, we set $\xi(x)=-V_{\omega}(x)$. 
In order to prove Theorem~\ref{annealed} by applying Theorem~\ref{GK},  
we first need to show that $H(t)$, which is obviously finite for all $t$, 
is close to the first term of the asymptotics. 
The following lemma establishes this, and also verifies Assumption (H). 
\begin{lemma}\label{lem1}
\begin{equation*}
	H(t) = -a_1 t^{\frac{d}{\alpha}}+O(e^{-t})
	 \quad {\rm as }\quad  t \to \infty.
\end{equation*}
\end{lemma}
\emph{Proof} 
It follows from \eqref{Laplace} that
\begin{equation*}
\begin{split}
	H(t) &= \log \E_{\nu} [\exp\{ -t V_{\omega}(0) \}] \\
	&= -\nu\int_{|y| \le 1}( 1-e^{-t}) {\rm d} y 
	- \nu\int_{|y| > 1} ( 1-e^{-t|y|^{-\alpha}} ) {\rm d} y\\ 
	&= O(e^{-t}) - \nu\int ( 1-e^{-t|y|^{-\alpha}} ) {\rm d} y. 
\end{split}
\end{equation*}
A computation shows that the second term in the 
last line equals $a_1 t^{\frac{d}{\alpha}}$.
\qed

\vspace{10pt} 
By this lemma, Theorem~\ref{annealed} turns out to be equivalent to 
\begin{equation*}
\begin{split}
	&\E_{\nu} \left[ E_0 \left[ \exp\set{-\int_0^t V_{\omega}(B_s) {\rm d} s} \right]^p \right]\\
	&\quad=\exp\set{ H(pt) - (a_2+o(1))(pt)^{\frac{\alpha+d-2}{2\alpha}}}. 
\end{split}
\end{equation*}

Our next task is to verify Assumption (J) for the scale
\begin{equation*}
	r=r(t)=t^{\frac{\alpha-d+2}{4\alpha}}
\end{equation*} 
and to identify the functional $J$. 
We first use \eqref{Laplace} to see
\begin{equation*}
\begin{split}
	&\log \E_{\nu}\left[ \exp\set{-t \int V_{\omega}(rx) \,\mu({\rm d} x)}\right] - H(t)\\
	&\quad = -\nu r^d \int \left(1- \exp\set{-t \int \hat{v}_r(x-y) \, \mu({\rm d} x)}\right){\rm d} y\\ 
	&\qquad\qquad +\nu r^d \int \left( 1-e^{-t \hat{v}_r(-y)} \right){\rm d} y,
\end{split}
\end{equation*}
where $\hat{v}_r(z)=\hat{v}(rz)$. 
Let us define $m_{\mu} = \int x \,\mu({\rm d} x)$ for $\mu \in \mathcal{P}_c(\R^d)$ 
and change the variable $y$ to $y-m_{\mu}$ in the second integral. 
Then, since both integrands are nonnegative and 
$ 1-e^{-t \hat{v}_r(m_{\mu}-y)}$ is integrable, we can merge the two integrals 
with respect to $y$ and arrive at 
\begin{equation*}
\begin{split}
	J_t(\mu) = \nu t^{-\frac{\alpha+d-2}{2\alpha}} r^d \int \left(e^{-t \hat{v}_r(m_{\mu}-y)}- 
	e^{-t \int \hat{v}_r(x-y) \, \mu({\rm d} x)} \right){\rm d} y. 
\end{split}
\end{equation*}
The following proposition verifies Assumption (J) and completes the proof of Theorem~\ref{annealed}.  
\begin{proposition}\label{prop1}
For any compact set $K \subset \R^d$, 
\begin{equation*}
	J_t(\mu) \xrightarrow{t \to \infty} \int \frac{\nu\alpha\sigma_d}{2d}  
	\Gamma\left( \frac{2\alpha-d+2}{\alpha} \right) |x-m_{\mu}|^2 \,\mu ({\rm d} x) 
\end{equation*}
uniformly in $\mu \in \mathcal{P}(K)$. 
 
\end{proposition}
\emph{Proof} 
We may assume $\nu=1$ and $m_{\mu}=0$ without loss of generality. 
Fix $R > 0$ such that $K \subset B(0,R)$. 
We first prove that relatively small $y$'s make only negligible contributions. 
\begin{lemma}\label{lem2}
For any $R > 0$, 
\begin{equation*}
	t^{-\frac{\alpha+d-2}{2\alpha}} r^d\int_{|y| \le t^{\frac{d+2-\alpha}{4(\alpha+1)}}} 
	\left(e^{-t \hat{v}_r(y)}- 
	e^{-t \int \hat{v}_r(x-y) \, \mu ({\rm d} x)} \right){\rm d} y \rightarrow 0 
\end{equation*}
as $t \to \infty$ uniformly in $\mu \in \mathcal{P}(B(0,R))$. 
\end{lemma}
\emph{Proof} 
The absolute value of the integrand is bounded from above by 
\begin{equation*}
\begin{split}
	& \max \set{e^{-t\hat{v}_r(y)}+e^{-t \int \hat{v}_r(x-y) \, \mu ({\rm d} x)}
	: \mu \in \mathcal{P}( B(0,R)), |y| \le t^{\frac{d+2-\alpha}{4(\alpha+1)}}}\\
	&\quad \le \exp\set{-tr^{-\alpha}\bigl(t^{\frac{d+2-\alpha}{4(\alpha+1)}}+R\bigr)^{-\alpha}}\\
	&\quad = \exp\set{-t^{\frac{d+2-\alpha}{4(\alpha+1)}}(1+o(1))}
\end{split}
\end{equation*}
as $t$ goes to infinity. Since $r^d$ and the volume of the integration range are both 
polynomial in $t$, the claim follows. 
\qed

\vspace{10pt}
This lemma allows us to consider only $y$'s with large modulus when $t$ is large. 
In what follows, we shall assume $t$ sufficiently large depending only on $\alpha$, $d$, and $R$, 
as necessary.  
Then $\hat{v}$ may be replaced by $v$ (recall $v(x)=|x|^{-\alpha}$) 
and we are reduced to proving 
\begin{equation*}
	t^{-\frac{\alpha+d-2}{2\alpha}} r^d\int_{|y| > t^{\frac{d+2-\alpha}{4(\alpha+1)}}} 
	\left(e^{-t v_r(y)}- e^{-t \int v_r(x-y) \, \mu ({\rm d} x)} \right){\rm d} y 
	\xrightarrow{t \to \infty} J (\mu)
\end{equation*}
uniformly in $\mu \in \mathcal{P}(B(0,R))$. The change of the variable $y=t^{1/\alpha}r^{-1}\eta$ 
shows that the above left hand side equals 
\begin{equation}
	t^{\frac{d+2-\alpha}{2\alpha}}
	\int_{|\eta| > t^{-\frac{d+2-\alpha}{4\alpha(\alpha+1)}}}e^{-v(\eta)}\left( 1- 
	e^{-\int (v(\eta-t^{-1/\alpha}rx)-v(\eta)) \, \mu ({\rm d} x)} \right){\rm d} \eta. 
	\label{scaled-prop1} 
\end{equation}
We use Taylor's theorem to approximate the integrand of $\int \mu({\rm d} x)$ as follows: 
\begin{equation*}
\begin{split}
	&v( \eta-t^{-1/\alpha}rx)-v(\eta) \\
	&\quad = -t^{-1/\alpha} r \langle \nabla v(\eta), x \rangle 
	+\frac{1}{2}t^{-2/\alpha} r^2\langle x, {\rm Hess}_v(\eta)x \rangle \\
	&\qquad \qquad 
	+ \int_0^1 \frac{(1-\theta)^2}{2} \frac{{\rm d}^3}{{\rm d}\theta^3}
	v(\theta t^{-1/\alpha}r x-\eta)\, {\rm d}\theta\\
	&\quad =: R_1(\eta, x)+ R_2(\eta, x) + R_3(\eta, x).  
\end{split}
\end{equation*}
It follows that $\int R_1(\eta, x)\,\mu({\rm d} x)=0$ from our assumption $m_{\mu}=0$. 
Moreover, considering the ranges of variables $x$ and $\eta$, one can easily see that 
\begin{equation}
	|R_i(\eta, x)| \le c_1(d,\alpha, R) |\eta|^{-\alpha-i} t^{-\frac{i}{\alpha}} r^i 
	\xrightarrow{t \to \infty} 0 
	\label{reminder-prop1} 
\end{equation}
for $i = 2,3$, where $c_1(d,\alpha, R)>0$ is a constant. 
In particular, the $\mu({\rm d} x)$ integral in \eqref{scaled-prop1} goes to 0 as $t \to\infty$. 
Hence we can use an elementary inequality 
$| 1-e^{-z}-z+z^2/2 | < |z|^3$
which holds when $|z|$ is small to obtain 
\begin{equation*}
\begin{split}
	&1 - e^{-\int (v(\eta-rt^{-1/\alpha}x)-v(\eta)) \, \mu ({\rm d} x)}\\
	&\quad =\int R_2(\eta, x) \, \mu ({\rm d} x)
	-\frac{1}{2}\left(\int R_2(\eta, x) \, \mu ({\rm d} x)\right)^2
	+ R_4(\eta), 
\end{split}
\end{equation*}
where $|R_4(\eta)| \le c_1^3t^{-\frac{3}{\alpha}} r^3 ( |\eta|^{-\alpha-3} \vee |\eta|^{-3\alpha-9})$. 

Now we perform the integration with respect to $e^{-v(\eta)}{\rm d} \eta$. 
First, Fubini's theorem and a little calculus show that 
\begin{equation}
\begin{split}
	&t^{\frac{d+2-\alpha}{2\alpha}} \int_{|\eta| > t^{-\frac{d+2-\alpha}{4\alpha(\alpha+1)}}}
	\int R_2(\eta, x) \, \mu ({\rm d} x)\, e^{-v(\eta)}{\rm d} \eta\\
	& \quad = \frac{1}{2}\int \left \langle x, \int_{|\eta| > t^{-\frac{d+2-\alpha}{4\alpha(\alpha+1)}}}
	{\rm Hess}_v(\eta) \, e^{-v(\eta)} {\rm d} \eta\, x \right \rangle\, \mu ({\rm d} x)\\
	& \quad \xrightarrow{t \to \infty} 
	\frac{\alpha\sigma_d}{2d}\Gamma\left( \frac{2\alpha-d+2}{\alpha} \right) 
	\int |x|^2 \,\mu ({\rm d} x)\label{a_2}
\end{split}
\end{equation}
uniformly in $\mu \in \mathcal{P}(B(0,R))$. 
Next, using the first inequality in \eqref{reminder-prop1} for $i=2$, we get 
\begin{equation*}
\begin{split}
	&t^{\frac{d+2-\alpha}{2\alpha}} \int_{|\eta| > t^{-\frac{d+2-\alpha}{4\alpha(\alpha+1)}}}
	\left(\int R_2(\eta, x) \, \mu ({\rm d} x)\right)^2 \, e^{-v(\eta)}{\rm d} \eta\\
	& \quad \le c_1^2 
	t^{-\frac{d+2-\alpha}{2\alpha}} \int |\eta|^{-2\alpha-4}e^{-|\eta|^{-\alpha}}{\rm d} \eta\\
	& \quad \xrightarrow{t \to \infty} 0
\end{split}
\end{equation*}
uniformly in $\mu \in \mathcal{P}(B(0,R))$. Finally, the integral of $R_4(\eta)$ can also be 
shown to converge to 0 uniformly by a similar estimate. 
\qed

\section{Upper bound on the integrated density of states}
\label{Lifshiz-upper}
We derive the upper bound on the integrated density of states. 
To this end, we employ the following well known relation (see e.g.~\cite{CL90}, Theorem VI.1.1): 
\begin{equation*}
\begin{split}
	\int_0^{\infty} e^{-t l} {\rm d} N(l) = (4\pi\kappa t)^{-\frac{d}{2}}\E_{\nu} \otimes E_{0,0}^t 
	\left[ \exp\set{-\int_0^t V_{\omega}(B_s) {\rm d} s} \right], 
\end{split}
\end{equation*} 
where $E_{0,0}^t$ denotes the expectation with respect to the Brownian bridge 
from $0$ to $0$ in the duration $t$. 
This right hand side is quite similar to the left hand side of \eqref{thm1}
and indeed exhibits the same asymptotic behavior. 
We need (and prove) only the following upper bound but the other direction can also be 
proved by the same argument as for Lemma~5 in~\cite{FU09}. 

\begin{lemma}\label{bridge}
\begin{equation*}
\begin{split}
	&(4\pi\kappa t)^{-\frac{d}{2}}\E_{\nu} \otimes E_{0,0}^t 
	\left[ \exp\set{-\int_0^t V_{\omega}(B_s) {\rm d} s} \right]\\
	&\quad \le \exp\set{H(t) -(a_2+o(1))t^{\frac{\alpha+d-2}{2\alpha}}}
\end{split}
\end{equation*}
as $t \to \infty$. 
\end{lemma}
\emph{Proof} 
By using a defining property of the Brownian bridge (see p.137 in~\cite{Szn98}), 
we find that the above left hand side 
is less than or equal to 
\begin{equation*}
	\E_{\nu} \otimes E_0
	\left[ \exp\set{-\int_0^{t-1} V_{\omega}(B_s) {\rm d} s} p(1, B_{t-1},0) \right], 
\end{equation*}
where 
\begin{equation*}
	p(t,x,y)=\frac{1}{(4\pi \kappa t)^{d/2}}\exp\set{-\frac{|x-y|^2}{4\kappa t}}. 
\end{equation*}
is the transition kernel of our Brownian motion. 
Since $p(1, B_{t-1}, 0) \le (4\pi \kappa )^{-d/2}$ and 
\begin{equation*}
	H(t-1)-a_2(t-1)^{\frac{\alpha+d-2}{2\alpha}}
	= H(t) - a_2t^{\frac{\alpha+d-2}{2\alpha}}+o(1)
\end{equation*}
as $t \to \infty$ by Lemma~\ref{lem1}, our claim follows from Theorem~\ref{annealed}. 
\qed

\vspace{10pt}
Using this lemma, we obtain 
\begin{equation}
\begin{split}
	N(\lambda) &\le e^{\lambda t} \int_0^{\lambda} e^{-t l} {\rm d} N(l) 
	\le e^{\lambda t} \int_0^{\infty} e^{-t l} {\rm d} N(l)\\
	&\le \exp\set{\lambda t -a_1t^{\frac{d}{\alpha}} -(a_2+o(1))t^{\frac{\alpha+d-2}{2\alpha}}}.
	\label{exact}
\end{split}
\end{equation}
Minimizing the first two terms of the right hand side over $t$, that is attained at
\begin{equation}
	\rho(\lambda)
	= \biggl( \frac{\alpha \lambda}{d a_1}\biggr)^{-\frac{\alpha}{\alpha-d}}, 
	\label{minimizer}
\end{equation}
we get the upper bound of Theorem~\ref{Lifshitz}. 

\section{Almost sure asymptotics}\label{proof-quenched}
In this section, we prove Theorem~\ref{quenched}. We first deal with the upper bound, which is rather easy, 
and then turn to more involved lower bound. \\

\noindent
\emph{Proof of the upper bound of Theorem~\ref{quenched}} 
The idea of the proof of the upper bound is close to that in~\cite{Fuk09b}.
We first show that the asymptotics of the Feynman--Kac functional can be controlled by 
the smallest Dirichlet eigenvalue of $-\kappa\Delta+V_{\omega}$ in a large box; 
we then derive an almost sure lower bound on the principal eigenvalue from 
Theorem~\ref{Lifshitz} by using a certain functional analytic inequality.  

We begin with the following general upper bound. 
\begin{lemma}\label{lem4}
There exist constants $c_2(d, \kappa), c_3(d, \kappa) > 0$ such that 
\begin{equation}
\begin{split}
	&E_0\left[ \exp\set{-\int_0^t V_{\omega}(B_s) {\rm d} s} \right]\\ 
	&\quad \le c_2 (1+(\lambda^{\rm D}_{\omega,\, 1}\bigl((-t, t)^d\bigr) t)^{d/2}) 
	\exp\set{ - \lambda^{\rm D}_{\omega,\, 1}\bigl((-t, t)^d\bigr) t} + e^{-c_3 t}. \label{spec}
\end{split}
\end{equation}
\end{lemma}
\emph{Proof} 
Let $\tau$ denote the exit time of the process from $(-t, t)^d$. 
Then, by using the reflection principle, one can show that there exists a constant
$c_3(d, \kappa)>0)$ such that
\begin{equation*}
\begin{split}
	 &E_0\left[ \exp\set{-\int_0^t V_{\omega}(B_s) {\rm d} s} \right]\\ 
	&\quad \le E_0 \left[\exp\set{-\int_0^t V_{\omega}(B_s){\rm d} s}: \tau > t\right]
	+ P_0 (\tau \le t)\\
	&\quad \le E_0 \left[\exp\set{-\int_0^t V_{\omega}(B_s){\rm d} s}: \tau > t\right]
	+ e^{-c_3 t}. 
\end{split}
\end{equation*}
Now, \eqref{spec} follows immediately from~(3.1.9) in p.93 of~\cite{Szn98}. 
\qed

\vspace{10pt}
Due to this lemma, it suffices to obtain the almost sure lower 
bound for the smallest Dirichlet eigenvalue $\lambda^{\rm D}_{\omega,\, 1}\bigl((-t, t)^d\bigr)$. 
We know from the fact mentioned in Remark~\ref{Rem-IDS}-(iii) that for any $\lambda > 0$ and $R > 0$,
\begin{equation}
\begin{split}
	N(\lambda) & \ge \frac{1}{(2R)^d} 
	\E_{\nu} \bigl[\#\bigl\{ k \in \N ;\lambda^{\rm D}_{\omega,\,k}\bigl((-R, R)^d\bigr) \le \lambda \bigr\}\bigr]\\
	& \ge \frac{1}{(2R)^d} \P_{\nu}\bigl(\lambda^{\rm D}_{\omega,\, 1}\bigl((-R, R)^d\bigr) 
	\le \lambda\bigr). \label{IDS-lower}
\end{split}
\end{equation} 
Now, let us fix $\e > 0$ arbitrarily and take 
\begin{equation}
	\lambda = \lambda(t,\e) = q_1 (\log t)^{-\frac{\alpha-d}{d}}
	+ \left(q_2 -\e \right)(\log t)^{-\frac{\alpha-d+2}{2d}}\label{lambda}
\end{equation}
and $R = t$. 
It then follows straightforwardly from Theorem~\ref{Lifshitz} and \eqref{IDS-lower} that 
\begin{equation*}
	\P_{\nu} \left( \lambda^{\rm D}_{\omega,\, 1}\bigl((-t, t)^d\bigr) \le \lambda(t,\e) \right) 
	\le \exp\set{-c_4(d,\nu,\alpha,\e)(\log t)^{\frac{\alpha+d-2}{2\alpha}}}
\end{equation*}
for some $c_4(d,\nu,\alpha,\e) > 0$ when $t$ is sufficiently large. 
This right hand side is summable along the sequence $t_k = e^k$ and hence Borel-Cantelli's 
lemma shows that $\P_{\nu}$-almost surely,
\begin{equation*}
	\lambda^{\rm D}_{\omega,\, 1}\bigl((-t_k, t_k)^d\bigr) \ge \lambda(t_k,\e)
\end{equation*} 
except for finitely many $k$. 
We can extend this bound for all large $t$ as follows: for $t_{k-1} \le t \le t_k$, we have 
\begin{equation*}
 \begin{split}
	\lambda^{\rm D}_{\omega,\, 1}\bigl((-t, t)^d\bigr) 
	& \ge \lambda^{\rm D}_{\omega,\, 1}\bigl((-t_k, t_k)^d\bigr) \\
	& \ge \lambda(t_k,\e) 
\end{split}
\end{equation*} 
by monotonicity and since 
\begin{align*}
	&(\log t_{k+1})^{-\frac{\alpha-d}{d}}-(\log t_k)^{-\frac{\alpha-d}{d}} 
	=o\bigl( (\log t_k)^{-\frac{\alpha-d+2}{2d}} \bigr), \\ 
	&(\log t_{k+1})^{-\frac{\alpha-d+2}{2d}} = (\log t_k)^{-\frac{\alpha-d+2}{2d}} (1+ o(1)) 
\end{align*}
as $t \to \infty$, we have $\lambda(t_k,\e) \ge \lambda(t,2\e)$ for sufficiently large $t$. 
Combined with Lemma~\ref{lem4}, this proves the upper bound of Theorem~\ref{quenched}. 
\qed
 
\vspace{10pt}\noindent
\emph{Proof of the lower bound of Theorem~\ref{quenched}} 
The proof of the lower bound of Theorem~\ref{quenched} goes along in the same spirit 
as~\cite{GKM00} and~\cite{Szn93c}. 
Roughly speaking, if we can find a sufficiently large \emph{pocket} not too far from the origin 
in which the smallest Dirichlet eigenvalue is close to 
$q_1(\log t)^{-(\alpha-d)/d} + q_2(\log t)^{-(\alpha-d+2)/2d}$, 
we can derive the lower bound 
by making the Brownian motion reach the pocket in relatively short time and stay there 
afterward. 
The following proposition gives the precise formulation and the almost sure existence of the pocket. 
\begin{proposition}\label{local}
For any $\e>0$ there exists $M>0$ such that the following holds: 
$\P_{\nu}$-almost surely, there exists a ball 
\begin{equation*}
	B_{\e,M}(t,\omega)=B \left(x_{\e,M}(t,\omega), M(\log t)^{\frac{\alpha-d+2}{4d}}\right), \;
	| x_{\e,M}(t,\omega) | \le t(\log t)^{-6}
\end{equation*}
for all sufficiently large $t$ such that
\begin{equation}
	\lambda^{\rm D}_{\omega, 1}(B_{\e,M}(t,\omega)) 
	\le q_1(\log t)^{-\frac{\alpha-d}{d}} + (q_2+\e)(\log t)^{-\frac{\alpha-d+2}{2d}}\label{pocket}.
\end{equation}
\end{proposition}
The proof of this proposition is slightly involved. 
We postpone it to the following subsections and first see 
how to derive the lower bound of Theorem~\ref{quenched}
from this. 
We need the following lemma to bound the potential on the way to the pocket. 
\begin{lemma}\label{potential-sup}
$\P_{\nu}$-almost surely, 
\begin{equation*}
	\sup_{x \in (-t, t)^d}V_{\omega}(x) \le  3d \log t 
\end{equation*} 
for sufficiently large $t$. 
\end{lemma}
\emph{Proof} 
Let $|\cdot|_{\infty}$ denote the $\ell^{\infty}$-norm on $\R^d$. 
We introduce the function $\bar{v}(x)=\sup_{| x-y |_{\infty} \le 1}\hat{v}(y)$ and 
\begin{equation*}
	\bar{V}_{\omega}(x)=\sum_i \bar{v}(x-\omega_i).
\end{equation*}
Then it is easy to see that 
\begin{equation*}
\begin{split}
	\E_{\nu}\biggl[\exp\biggl\{ \sup_{x \in [0,1)^d}V_{\omega}(x)\biggr\}\biggr]
	& \le \E_{\nu}[\exp\{ \bar{V}_{\omega}(0)\}]\\
	& = \exp\set{\nu\int (e^{\bar{v}(-y)}-1){\rm d} y}\\
	& < \infty. 
\end{split}
\end{equation*}
Therefore, Chebyshev's inequality shows 
\begin{equation*}
\begin{split}
	\P_{\nu} & \biggl(\sup_{x \in (-2t,2t)^d} V_{\omega}(x) > 3d \log t \biggr)\\
	& \le (4t)^d \P_{\nu} \biggl(\sup_{x \in [0,1)^d} V_{\omega}(x) > 3d \log t \biggr)\\
	& \le 4^d t^{-2d}\, \E_{\nu}\biggl[\sup_{x \in [0,1)^d}\exp\{V_{\omega}(x)\}\biggr].
\end{split}
\end{equation*}
Since the last expression is summable in $t \in \N$, the claim follows by Borel-Cantelli's lemma 
and monotonicity of $\sup_{x \in (-t, t)^d}V_{\omega}(x)$ in $t$. 
\qed

\vspace{10pt}
Now let us pick $\omega$ and sufficiently large $t$ so that assertions in Proposition~\ref{local} and 
Lemma~\ref{potential-sup} holds. 
We denote by $\phi_{\omega}$ the $L^2$-normalized nonnegative eigenfunction 
associated with $\lambda^{\rm D}_{\omega, 1}(B_{\e,M}(t,\omega))$. 
Then, since we know the following uniform upper bound from (3.1.55) in~\cite{Szn98} 
\begin{equation*}
	\| \phi_{\omega}\|_{\infty} \le c_5(d) \lambda^{\rm D}_{\omega,\, 1}(B_{\e,M}(t,\omega))^{d/4}, 
\end{equation*}
the integral of $\phi_{\omega}$ is bounded from below as 
\begin{equation*}
	\int \phi_{\omega}(x) {\rm d} x 
	\ge \frac{1}{\| \phi_{\omega}\|_{\infty} }\int \phi_{\omega}(x)^2 {\rm d} x 
	\ge c_5^{-1}q_1^{-\frac{d}{4}}(\log t)^{-\frac{\alpha-d}{4}}.
\end{equation*}
Now, recall that the Feynman-Kac semigroup generated by $-\kappa \Delta + V_{\omega}$ has the 
kernel $p_{\omega}(s,x,y)$ since the potential term is locally bounded (see Theorem B.7.1 in~\cite{Sim82}). 
We can bound this kernel from below by using the Dirichlet heat kernel $p_{(-t, t)^d}(s,x,y)$ 
in $(-t, t)^d$ as follows: 
\begin{equation*}
\begin{split}
	p_{\omega}(s,0,y) 
	&\ge \exp\Bigl\{-s \sup_{x \in (-t, t)^d}V_{\omega}(x)\Bigr\} p_{(-t, t)^d}(s,0,y)\\
	& \ge c_6 s^{-d/2}\exp\bigl\{-s (3d \log t) - c_6^{-1}{|y|^2}/{s} \bigr\} 
	\quad \textrm{if} \quad |y|<t/2,  
\end{split}
\end{equation*}
where $c_6(d, \kappa)$ is a constant and 
the second inequality follows by Lemma~\ref{potential-sup} and a Gaussian lower bound for 
the Dirichlet heat kernel in~\cite{vdB92}. 
Taking $s=t(\log t)^{-6}$ and noting that $B_{\e,M}(t,\omega) \subset 
(-2t(\log t)^{-6}, 2t(\log t)^{-6})^d$, we obtain 
\begin{equation}
	\inf_{y \in B_{\e,M}(t,\omega)} p_{\omega}(t(\log t)^{-6},0,y) 
	\ge \exp\set{ -4d t(\log t)^{-5} } \label{kernel}
\end{equation}
for sufficiently large $t$. 

By Chapman-Kolmogorov's equation, we have
\begin{equation*}
\begin{split}
	 &E_0\left[ \exp\set{-\int_0^t V_{\omega}(B_s) {\rm d} s} \right]
	 = \int p_{\omega}(t,0,x) {\rm d} x\\
	& \quad\ge \iint p_{\omega}(t(\log t)^{-6}, 0, y) p_{\omega}(t-t(\log t)^{-6}, y, x)  
	\frac{\phi_{\omega}(x)}{\|\phi_{\omega} \|_{\infty}}{\rm d} y {\rm d} x.
\end{split}
\end{equation*}
We use~ \eqref{kernel} for the first $p_{\omega}$ in the second line and replace the second $p_{\omega}$ by 
the kernel of the semigroup generated by $-\kappa \Delta+V_{\omega}$ with the Dirichlet boundary condition 
outside $B_{\e,M}(t,\omega)$. 
Then, since $\phi_{\omega}$ is the eigenfunction and we have \eqref{pocket}, 
we find that the above right hand side is bounded from below by
\begin{equation*}
\begin{split}
	& \exp\bigl\{ -\lambda_{\omega,\, 1}(B_{\e,M}(t,\omega))t - 5d t(\log t)^{-5} \bigr\}
	\frac{1}{\|\phi_{\omega}\|_{\infty}} \int \phi_{\omega}(y) {\rm d} y\\
	& \quad \ge \exp\set{ -q_1 t(\log t)^{-\frac{\alpha-d}{d}} - (q_2+2\e) t(\log t)^{-\frac{\alpha-d+2}{2d}}} 
\end{split}
\end{equation*}
for sufficiently large $t$. This completes the proof of the lower bound of Theorem~\ref{quenched}.
\qed
\subsection{Proof of Proposition~\ref{local} 
in the case $\alpha \ge 2$}
We prove Proposition~\ref{local} in the case $\alpha \ge 2$ in this
subsection. 
Since the remaining case $1<\alpha<2$ is treated in a similar way, 
we specify where we need $\alpha \ge 2$ and only give necessary changes
in Subsection~\ref{remaining}.

We first introduce a quadratic function 
\begin{equation*}
	Q_t(x)= q_1(\log t)^{-\frac{\alpha-d}{d}}+
	\frac{q_2^2}{\kappa d}(\log t)^{-\frac{\alpha-d+2}{d}}|x|^2. 
\end{equation*}
and a ball
\begin{equation*}
	B_M(t)=B\left(0, M(\log t)^{\frac{\alpha-d+2}{4d}}\right). 
\end{equation*}
Then, it is not difficult to see by using a scaling that the smallest Dirichlet eigenvalue of $-\kappa \Delta + Q_t$ 
in $B_M(t)$ is 
\begin{equation*}
	q_1 (\log t)^{-\frac{\alpha-d}{d}}
	+ q_2 (\log t)^{-\frac{\alpha-d+2}{2d}}(1+o(1)) 
\end{equation*}
as $M \to \infty$, uniformly in $t$. 
Hence, if we show that there exists a pocket $B_{\e,M}(t,\omega)$ with large $M$ in which $V_{\omega}$
is close to a translation of $Q_t(x)$, Proposition~\ref{local} follows. 
This is indeed possible in the case $\alpha \ge 2$. 
\begin{proposition}\label{alpha>2}
Suppose $\alpha \ge 2$. Then for any $M>0$ and $\e>0$, there $\P_{\nu}$-almost surely exists a ball 
\begin{equation}
	B_{\e,M}(t,\omega)=x_{\e,M}(t,\omega)+B_M(t), \;
	| x_{\e,M}(t,\omega) | \le t(\log t)^{-6}\label{ball2}
\end{equation}
for all sufficiently large $t$ in which we have 
\begin{equation}
	|V_{\omega}(x) - Q_t(x-x_{\e,M}(t,\omega))| \le \e(\log t)^{-\frac{\alpha-d+2}{2d}}.
	\label{approx-potential}
\end{equation}
\end{proposition}
\emph{Proof} 
Let us first explain the strategy of the proof. We first choose a collection of balls 
which are 
\begin{itemize}
\item{translations of $B_M(t)$ like~\eqref{ball2}, }
\item{more than $t^d(\log t)^{-C}$ in number for some $C>0$,}
\item{so far from each other that the shapes of $V_{\omega}$ 
in different balls are almost independent. }
\end{itemize}
Next, we show that in each ball, the probability of \eqref{approx-potential} 
is larger than $t^{-d}\exp\{(\log t)^{\delta}\}$ for some $\delta \in (0,1)$. 
Then, roughly speaking, we have $t^d(\log t)^{-C}$ almost independent trials with 
success probability more than $t^{-d}\exp\{(\log t)^{\delta}\}$ and this assures at least one success. 

Now we go into the rigorous argument. We fix  
\begin{equation*}
	N > \frac{\alpha-d+2}{2d(\alpha-d)}\vee 2
\end{equation*}
and then set 
\begin{equation*}
	\I = 2(\log t)^N \Z^d,
\end{equation*}
which will be the centers of the collection of balls. 
We write $\Lambda_N$ for $[-(\log t)^N,(\log t)^N)^d$ to simplify notation. 
The following lemma corresponds to the independence property in the 
first step of the strategy. 
\begin{lemma}\label{lem6}
For any $\e >0$, 
\begin{equation}
 \P_{\nu}\left(\sup_{y \in B_M(t)} 
	\sum_{\omega_i \not\in \Lambda_N}|y-\omega_i|^{-\alpha}
	> \e (\log t)^{-\frac{\alpha-d+2}{2d}} \right)
 	\le \exp\set{-(\log t)^{dN}}.\label{indep}
\end{equation}
\end{lemma}
\noindent\emph{Proof} 
For $\omega_i \not\in \Lambda_N$ and $y \in B_M(t)$, 
we have $|y-\omega_i|>|\omega_i|/2$ and thus
\begin{equation*}
 \sup_{y \in B_M(t)}
 \sum_{\omega_i \not\in \Lambda_N}|y-\omega_i|^{-\alpha}
 < 2^{\alpha}\sum_{\omega_i \not\in \Lambda_N}
 \left|\omega_i\right|^{-\alpha}.
\end{equation*}
We use~\eqref{Laplace} to see 
\begin{equation*}
  \E_{\nu}\left[\exp\set{(\log t)^{\alpha N}
  \sum_{\omega_i \not\in \Lambda_N}
  \left|\omega_i\right|^{-\alpha}}\right]
  =\exp\set{\nu\int_{\R^d\setminus\Lambda_N}
  (e^{(\log t)^{\alpha N}|z|^{-\alpha}}-1)\d z}.
\end{equation*}
Now since $(\log t)^{\alpha N}|z|^{-\alpha}$ is bounded for 
$z\not\in \Lambda_N$, we have 
$e^{(\log t)^{\alpha N}|z|^{-\alpha}}-1 
\le c_7(\log t)^{\alpha N}|z|^{-\alpha}$ for 
some $c_7(d,\alpha)>0$ 
and the above integral is bounded as
\begin{equation*}
 \begin{split}
  \int_{\R^d\setminus\Lambda_N}
  (e^{(\log t)^{\alpha N}|z|^{-\alpha}}-1)\d z
  &\le c_7\int_{\R^d\setminus\Lambda_N}
  (\log t)^{\alpha N}|z|^{-\alpha}\d z\\
  &=O\left((\log t)^{(\alpha-d) N}\right).
 \end{split}
\end{equation*}
Then Chebyshev's inequality yields
\begin{equation*}
 \begin{split}
  &\P_{\nu}\left( 
	2^{\alpha}
        \sum_{\omega_i \not\in \Lambda_N}|y-\omega_i|^{-\alpha}
	> \e (\log t)^{-\frac{\alpha-d+2}{2d}} \right)\\
 &\le \exp\set{-2^{-\alpha}\epsilon(\log t)^{\alpha N
               -\frac{\alpha-d+2}{2d}}
               +O\left((\log t)^{(\alpha-d) N}\right)}
 \end{split}
\end{equation*}
and this implies~\eqref{indep} thanks to our choice of $N$.
\qed

\vspace{10pt}
To estimate the probability of each ball being a pocket, it is 
convenient to introduce 
a transformed measure defined by 
\begin{equation*}
	\frac{{\rm d}\Tilde{\P}_{t}}{{\rm d}\P_{\nu}}(\omega) 
	= e^{-H(\rho(\lambda(t)))-\rho(\lambda(t))V_{\omega}(0)}, 
\end{equation*}
where $\rho$ is defined in \eqref{minimizer} and 
\begin{equation*}
	\lambda(t)=q_1(\log t)^{-\frac{\alpha-d}{d}}=Q_t(0). 
\end{equation*}
Note that with this choice
\begin{align}
	&\rho(\lambda(t))=\left(a_1\frac{\alpha-d}{\alpha d}\right)^{-\frac{\alpha}{d}}(\log t)^{\frac{\alpha}{d}},
	\label{rho(lambda)}\\
	&H(\rho(\lambda(t))+\lambda(t)\rho(\lambda(t))=-d\log t+o(1),\label{dlogt}
\end{align}
where the latter follows from Lemma~\ref{lem1}. 
Recall that $\rho$ defined in~\eqref{minimizer} is the minimizer of~\eqref{exact}.
We collect several properties of the measure~$\Tilde{\P}_{t}$ which we shall use later. 
\begin{lemma}\label{lem7}
\begin{enumerate}
\item{$(\omega, \Tilde{\P}_t)$ is a Poisson point process with intensity 
$\nu e^{-\rho(\lambda(t))\hat{v}(y)}{\rm d} y$. \\}
\item{$\Tilde{\E}_t[V_{\omega}(x)]=Q_t(x)+o((\log t)^{-\frac{\alpha-d+2}{2d}})$ as $t \to \infty$, 
uniformly in $x \in B_M(t)$.\\}
\item{ 
$(\log t)^{\frac{2\alpha-d}{2d}}(V_{\omega}(0) - \lambda(t))$ under $\Tilde{\P}_t$ 
converges in law to a non-degenerate Gaussian random variable.}
\end{enumerate}
\end{lemma}
\emph{Proof} 
The proof of (i) is straightforward. Indeed, for any nonnegative Borel function $f$ on $\R^d$, 
we have 
\begin{equation*}
\begin{split}
	& \Tilde{\E}_t\left[ \exp\set{-\int f(y) \,\omega({\rm d} y)} \right]\\
	& \quad = e^{-H(\rho(\lambda(t)))}
	\E_{\nu}\left[ \exp\set{\int (-f(y)-\rho(\lambda(t))\hat{v}(y)) \,\omega({\rm d} y)} \right]\\
	& \quad = \exp\set{-\nu \int (1-e^{-f(y)}) e^{-\rho(\lambda(t))\hat{v}(y)} {\rm d} y}
\end{split}
\end{equation*}
by using \eqref{Laplace}. This verifies a condition to identify a point process and (i) follows 
(see, e.g.,\ Proposition~3.6 of~\cite{Res87}). 

To prove (ii), note first that by \eqref{expectation} we have 
\begin{equation}
\begin{split}
	\Tilde{\E}_t[V_{\omega}(x)] &=\nu \int \hat{v}(x-y)e^{-\rho(\lambda(t))\hat{v}(-y)}{\rm d} y\\
	&=\nu\int_{B_{2M}(t)} \hat{v}(x-y)e^{-\rho(\lambda(t))\hat{v}(-y)}{\rm d} y\\
	&\qquad+\nu\int_{\R^d \setminus B_{2M}(t)}v(x-y)e^{-\rho(\lambda(t))v(-y)}{\rm d} y.
	\label{mean}
\end{split}
\end{equation}
(Recall $v(x)=|x|^{-\alpha}$.)
We say that a function $f(t)$ is of order $o((\log t)^{-\infty})$ if 
\begin{equation*}
	f(t)= o((\log t)^{-L})  
\end{equation*}
for any $L>0$. 
We shall use the following lemma in the sequel. 
\begin{lemma}\label{lem8}
\begin{enumerate}
\item{For any $M>0$, 
\begin{equation*}
	\sup_{\|u\|_{\infty} \le 1 }\left|\int_{B_{2M}(t)} u(y) e^{-\rho(\lambda(t))\hat{v}(y)}{\rm d} y\right| 
	= o((\log t)^{-\infty})
\end{equation*}
as $t \to \infty$.}
\item{For any $M>0$ and $\gamma >0$, 
\begin{equation*}
	\int_{B_{2M}(t)} |y|^{-\gamma} e^{-\rho(\lambda(t))v(y)}{\rm d} y = o((\log t)^{-\infty})
\end{equation*}
as $t \to \infty$. }
\item{For any $M>0$ and $\gamma >d$, 
\begin{equation}
	\int_{\R^d \setminus B_{2M}(t)} |y|^{-\gamma} 
	e^{-\rho(\lambda(t))\hat{v}(y)}{\rm{d}}y = O\left((\log t)^{\frac{d-\gamma}{d}}\right)\label{Lem8-3}
\end{equation}
as $t \to \infty$. }
\end{enumerate}
\end{lemma}
We omit the proof of this lemma since it is elementary.  
By using Lemma~\ref{lem8}-(i) for the first term in the right hand side of \eqref{mean}, we get
\begin{equation}
	\Tilde{\E}_t[V_{\omega}(x)] 
	= \nu\int_{\R^d \setminus B_{2M}(t)}v(x-y)e^{-\rho(\lambda(t))v(-y)}{\rm d} y
	+o((\log t)^{-\infty})\label{inside-negligible}
\end{equation}
as $t \to \infty$, uniformly in $x \in B_M(t)$. 
Moreover, it also follows from Lemma~\ref{lem8}-(ii) that
\begin{equation}
	\int_{B_{2M}(t)} v(-y)e^{-\rho(\lambda(t))v(-y)}{\rm d} y =o((\log t)^{-\infty}).\label{inside-replace}
\end{equation}
Therefore, by adding \eqref{inside-replace} to \eqref{inside-negligible} with $x=0$,  we obtain 
\begin{equation}
\begin{split}
	\Tilde{\E}_t[V_{\omega}(0)]
	&=\nu \int v(-y)e^{-\rho(\lambda(t))v(-y)}{\rm d} y + o((\log t)^{-\infty})\\
	&= \lambda(t) + o((\log t)^{-\infty})\label{atorigin}
\end{split}
\end{equation}
and the claim is proved for $x=0$. 
Next, for general $x \in B_M(t)$, we have  
\begin{equation*}
\begin{split}
	&\Tilde{\E}_t[V_{\omega}(x)-V_{\omega}(0)]\\
	&\quad =\nu \int_{\R^d \setminus B_{2M}(t)} (v(x-y)-v(-y))e^{-\rho(\lambda(t))v(-y)}{\rm d} y 
	+o((\log t)^{-\infty})
\end{split}
\end{equation*}
as $t \to \infty$ by \eqref{inside-negligible}. 
Applying Taylor's theorem to the integrand, we obtain
\begin{equation*}
\begin{split}
	& v(x-y)-v(-y)\\
	&\quad =\langle\nabla v(-y), x\rangle 
	+\frac{1}{2}\langle x, {\rm Hess}_v(-y)x \rangle 
	+ \int_0^1 \frac{(1-\theta)^2}{2}\frac{{\rm d}^3}{{\rm d}\theta^3}v(\theta x-y){\rm d} \theta.
\end{split}
\end{equation*}
The first term in the right hand side vanishes integrated against $\nu e^{-\rho(\lambda(t))v(-y)}{\rm d} y$ over 
$\R^d \setminus B_{2M}(t)$ since 
$\nabla v(-y)e^{-\rho(\lambda(t))v(-y)}$ and the integral region 
are symmetric about the origin. For the second term, one can show that
\begin{equation*}
\begin{split}
	&\frac{\nu}{2}
	\int_{\R^d \setminus B_{2M}(t)} \langle x, {\rm Hess}_v(-y)x \rangle e^{-\rho(\lambda(t))v(-y)}{\rm d} y \\
	&\quad =\frac{\nu}{2}\int \langle x, {\rm Hess}_v(-y)x \rangle e^{-\rho(\lambda(t))v(-y)}{\rm d} y
	+o((\log t)^{-\infty}) \\
	&\quad= \frac{q_2^2}{\kappa d}(\log t)^{-\frac{\alpha-d+2}{d}}|x|^2 +o((\log t)^{-\infty})
\end{split}
\end{equation*}
as $t \to \infty$ by using Lemma~\ref{lem8}-(ii) and the same calculation as in \eqref{a_2}. 
Finally, for the third term, note that there exists a constant $c_8(d,\alpha)>0$ such that
\begin{equation}
	\left|\frac{{\rm d}^3}{{\rm d}\theta^3}v(\theta x-y)\right|
	\le c_8M^3(\log t)^{\frac{3(\alpha-d+2)}{4d}}|y|^{-\alpha-3}\label{3rd-derivative}
\end{equation}
for any $\theta \in [0,1]$, $x \in B_M(t)$, and $y \not\in B_{2M}(t)$. 
Therefore, we can bound its integral as
\begin{equation*}
\begin{split}
	&\left|\int_{\R^d \setminus B_{2M}(t)}
	\int_0^1 \frac{(1-\theta)^2}{2}\frac{{\rm d}^3}{{\rm d}\theta^3}v(\theta x-y){\rm d} \theta
	e^{-\rho(\lambda(t))v(-y)}{\rm d} y \right|\\
	&\quad \le  \int c_8 |y|^{-\alpha-3}e^{-\rho(\lambda(t))v(-y)}{\rm d} y\\
	&\quad =O((\log t)^{-\frac{\alpha-d+6}{4d}})
\end{split}
\end{equation*}
by using the change of variable $y=\rho(\lambda(t))^{1/\alpha}\eta$.
The last line is of order $o((\log t)^{-(\alpha-d+2)/2d})$ and the proof of (ii) is completed.

To prove (iii), we first replace $\lambda(t)$ in the statement by $\Tilde{\E}_t[V_{\omega}(0)]$. 
This cause no difference thanks to \eqref{atorigin}. 
Then we use the first line of \eqref{mean} and \eqref{characteristic} to obtain 
\begin{equation*}
\begin{split}
	& \Tilde{\E}_t\left[ \exp\set{i \theta \rho(\lambda(t))^{\frac{2\alpha-d}{2\alpha}}
	(V_{\omega}(0)-\Tilde{\E}_t[V_{\omega}(0)])} \right]\\
	&\quad = \exp\Biggl\{\nu 
	\int \Bigl(e^{i \theta\rho(\lambda(t))^{\frac{2\alpha-d}{2\alpha}}\hat{v}(y)}-1\\
	&\hspace{80pt}-i \theta \rho(\lambda(t))^{\frac{2\alpha-d}{2\alpha}}\hat{v}(y)\Bigr) 
	e^{-\rho(\lambda(t))\hat{v}(y)} {\rm d} y\Biggr\}.
\end{split}
\end{equation*}
Let us write
\begin{equation*}
	F(\theta, t, y)= e^{i \theta\rho(\lambda(t))^{\frac{2\alpha-d}{2\alpha}}\hat{v}(y)}-1
	-i \theta \rho(\lambda(t))^{\frac{2\alpha-d}{2\alpha}}\hat{v}(y)
\end{equation*}
to simplify the notation. Since $F(\theta, t, y)$ is bounded in $y$ for fixed $t$ and grows polynomially in $t$, 
it is easy to see that for any $\beta \in (0, 1)$, 
\begin{equation*}
	\int_{|y| \le \rho(\lambda(t))^{\frac{1-\beta}{\alpha}}} 
	|F(\theta, t, y)| 
	e^{-\rho(\lambda(t))\hat{v}(y)} {\rm d} y
	=O\left(\exp\set{-\rho(\lambda(t))^{\beta}/2}\right) 
\end{equation*}
as $t \to \infty$. 
Hence this region makes only negligible contribution to the integral. 
On the other hand, we may replace $\hat{v}(y)$ by $v(y)$ on $\{|y| >\rho(\lambda(t))^{(1-\beta)/\alpha}\}$ 
and then the change of variable $y=\rho(\lambda(t))^{1/\alpha}\eta$ yields 
\begin{equation}
\begin{split}
	&\int_{|y| >\rho(\lambda(t))^{\frac{1-\beta}{\alpha}}} 
	F(\theta, t, y) e^{-\rho(\lambda(t))\hat{v}(y)} {\rm d} y\\
	& \quad = \rho(\lambda(t))^{d/\alpha}
	\int_{|\eta| > \rho(\lambda(t))^{-\frac{\beta}{\alpha}}} 
	\Bigl( e^{i \theta\rho(\lambda(t))^{-\frac{d}{2\alpha}}v(\eta)}-1\\
	&\hspace{135pt} -i \theta \rho(\lambda(t))^{-\frac{d}{2\alpha}} v(\eta) \Bigr) 
	e^{-v(\eta)} {\rm d} \eta. \label{ch.f.-main}
\end{split}
\end{equation}
Now we take $\beta < d/(2\alpha)$ so that 
\begin{equation*}
	\rho(\lambda(t))^{-\frac{d}{2\alpha}}v(\eta) \xrightarrow{t \to \infty} 0
\end{equation*}
uniformly on $\{|\eta| > \rho(\lambda(t))^{-\beta/\alpha}\}$. 
Then, we have
\begin{equation*}
\begin{split}
	&e^{i \theta\rho(\lambda(t))^{-\frac{d}{2\alpha}}v(\eta)}-1
	-i \theta \rho(\lambda(t))^{-\frac{d}{2\alpha}} v(\eta)\\
	&\quad =-\frac{\theta^2}{2} \rho(\lambda(t))^{-\frac{d}{\alpha}} v(\eta)^2(1+o(1)),
\end{split}
\end{equation*}
where $o(1)$ is uniform in $\eta$. 
Substituting this into~\eqref{ch.f.-main} and recalling 
\eqref{rho(lambda)}, we complete the proof of (iii). 
\qed

\vspace{10pt}
The next lemma establishes the second step of the strategy. 
\begin{lemma}\label{lem9}
Suppose $\alpha \ge 2$. 
Then there exists $\delta \in (0,1)$ such that for any $\e > 0$ and $M>0$, 
\begin{equation}
\begin{split}
	\P_{\nu}\left( \sup_{x \in B_M(t)} |V_{\omega}(x) - Q_t(x)| \le \e(\log t)^{-\frac{\alpha-d+2}{2d}}\right)
	\ge t^{-d}\exp\set{ (\log t)^{\delta}}\label{prob-lower}
\end{split}
\end{equation}
when $t$ is sufficiently large. 
\end{lemma}
\emph{Proof}
We introduce a parameter 
\begin{equation}
	\gamma \in \left( \frac{\alpha-d+2}{2d}, \frac{\alpha}{d} \right). 
\end{equation}
Now, in view of Lemma~\ref{lem7}-(ii), we have an inclusion
\begin{equation*}
\begin{split}
	&\set{ \sup_{x \in B_M(t)}|V_{\omega}(x) - Q_t(x)| \le \e(\log t)^{-\frac{\alpha-d+2}{2d}}}\\
	&\quad \supset \set{V_{\omega}(0) - \lambda(t)\in 
	\left((\log t)^{-\gamma}, \frac{\e}{2}(\log t)^{-\frac{\alpha-d+2}{2d}}\right)}\setminus \\ 
	&\quad \qquad \set{\sup_{x \in B_M(t)}|V_{\omega}(x)-V_{\omega}(0)-
	\Tilde{\E}_t[V_{\omega}(x)-V_{\omega}(0)]|
	\ge \frac{\e}{4}(\log t)^{-\frac{\alpha-d+2}{2d}}}\\
	&\quad=: E_1 \setminus E_2 
\end{split}
\end{equation*}
for sufficiently large $t$. From this and~\eqref{dlogt}, it follows that 
\begin{equation*}
\begin{split}
	&\textrm{the left hand side of \eqref{prob-lower}} \\
	&\quad \ge e^{H(\rho(\lambda(t)))}\Tilde{\E}_t \bigl[e^{\rho(\lambda(t))V_{\omega}(0)}: 
	E_1 \setminus E_2\bigr]\\
	&\quad \ge \exp\set{H(\rho(\lambda(t)))+\rho(\lambda(t))
	\bigl(\lambda(t) + (\log t)^{-\gamma}\bigr)}\Tilde{\P}_t(E_1\setminus E_2)\\
	&\quad \ge \exp\set{-d \log t + \rho(\lambda(t))(\log t)^{-\gamma}+o(1)}
	(\Tilde{\P}_t(E_1)-\Tilde{\P}_t(E_2)).
\end{split}
\end{equation*}
Since $\rho(\lambda(t))(\log t)^{-\gamma}$ is a positive power of $\log t$, 
it remains to show that $\Tilde{\P}_t(E_1)-\Tilde{\P}_t(E_2)$ is bounded from below. 
The first term is rather easy since 
\begin{equation}
\begin{split}
	\Tilde{\P}_t(E_1)=\Tilde{\P}_t \Bigl((\log t)^{\frac{2\alpha-d}{2d}}(V_{\omega}(0) - \lambda(t)) \in 
	\left((\log t)^{\frac{2\alpha-d}{2d}-\gamma}, \frac{\e}{2}(\log t)^{\frac{\alpha-2}{2d}}\right)\Bigr), 
	\label{E_1}
\end{split}
\end{equation}
which is bounded from below by a positive constant for $\alpha \ge 2$ because of Lemma~\ref{lem7}-(iii). 
To estimate $\Tilde{\P}_t(E_2)$, we use~\eqref{expectation} to see
\begin{equation}
\begin{split}
	&V_{\omega}(x)-V_{\omega}(0)-\Tilde{\E}_t[V_{\omega}(x)-V_{\omega}(0)]\\
	&\quad=\int (\hat{v}(x-y)-\hat{v}(-y)) 
	\left(\omega({\rm{d}}y) - \nu e^{-\rho(\lambda(t))\hat{v}(y)} {\rm{d}}y\right)\label{deviation}.
\end{split}
\end{equation}
For abbreviation, we write $\bar{\omega}_t({\rm{d}}y)$ for $\omega({\rm{d}}y)
-\nu e^{-\rho(\lambda(t))\hat{v}(y)} {\rm{d}}y$ in this proof. 
This is a slight abuse of notation since  $\bar{\omega}_t({\rm{d}}y)$ has 
infinite total variation. But we will only consider functions 
which are $\nu e^{-\rho(\lambda(t))\hat{v}(y)} {\rm{d}}y$-integrable and therefore all the 
integrals appearing below make sense. 

We divide the integral in~\eqref{deviation} into $y \in B_{2M}(t)$ and $y \not\in B_{2M}(t)$
and show that each part has order $o((\log t)^{-(\alpha-d+2)/2d})$ with probability 
close to 1. Fix an arbitrary small $\e>0$. 
Let us begin with
\begin{equation*}
\begin{split}
	&\sup_{x \in B_M(t)}\left|\int_{B_{2M}(t)} 
	(\hat{v}(x-y)-\hat{v}(-y)) \bar{\omega}_t({\rm{d}}y)\right|\\
	&\quad \le \sup_{x \in B_M(t)} \biggl\{ \int_{B_{2M}(t)}
	|\hat{v}(x-y)-\hat{v}(-y)| \omega({\rm{d}}y)\\
	&\hspace{80pt}  +\nu\int_{B_{2M}(t)} 
	|\hat{v}(x-y)-\hat{v}(-y)| e^{-\rho(\lambda(t))\hat{v}(y)} {\rm{d}}y\biggr\}\\
	&\quad \le  \int_{B_{2M}(t)}
	\bar{\omega}_t({\rm{d}}y)
	+2\nu\int_{B_{2M}(t)} e^{-\rho(\lambda(t))\hat{v}(y)} {\rm{d}}y.
\end{split}
\end{equation*}
The $\Tilde\P_t$-mean of the first term is zero. Moreover, its variance and 
the second term are both of $o((\log t)^{-\infty})$ by Lemma~\ref{lem8}-(i). 
Hence we obtain 
\begin{equation*}
\begin{split}
	&\Tilde\P_t\left(\sup_{x \in B_M(t)}\left|\int_{B_{2M}(t)} 
	(\hat{v}(x-y)-\hat{v}(-y)) \bar{\omega}_t({\rm{d}}y)\right| 
	> \e (\log t)^{-\frac{\alpha-d+2}{2d}}\right)\\
	&\quad=o((\log t)^{-\infty})
\end{split}
\end{equation*}
as $t \to \infty$ using Chebyshev's inequality. 

Now we turn to the remaining part. Since $\hat{v}(x-y)=v(x-y)(=|x-y|^{-\alpha})$ for 
$x \in B_M(t)$ and $y \not\in B_{2M}(t)$, we can use Taylor's theorem to see
\begin{equation}
\begin{split}
	&\sup_{x \in B_M(t)}\left|\int_{\R^d \setminus B_{2M}(t)}
	(\hat{v}(x-y)-\hat{v}(-y)) \bar{\omega}_t({\rm{d}}y)\right|\\
	&\quad = \sup_{x \in B_M(t)}\left|\int_{\R^d \setminus B_{2M}(t)}
	\langle x, \nabla v(-y)\rangle \bar{\omega}_t({\rm{d}}y)\right|\\
	&\qquad +\sup_{x \in B_M(t)}\left|\int_{\R^d \setminus B_{2M}(t)}
	\frac{1}{2}\langle x, {\rm Hess}_{v}(-y)x \rangle  
	\bar{\omega}_t({\rm{d}}y)\right|\\
	&\qquad +\sup_{x \in B_M(t)}\left|\int_{\R^d \setminus B_{2M}(t)}\int_0^1\frac{(1-\theta)^2}{2}
	\frac{{\rm d}^3}{{\rm d}\theta^3}v(\theta x-y) {\rm d}\theta
	\bar{\omega}_t({\rm{d}}y)\right|. \label{B_{2M}(t)^c}
\end{split}
\end{equation}
The first term on the right hand side is bounded as 
\begin{equation*}
\begin{split}
	&\sup_{x \in B_M(t)}\left|\int_{\R^d \setminus B_{2M}(t)}
	\langle x, \nabla v(-y)\rangle \bar{\omega}_t({\rm{d}}y)\right|\\
	&\quad \le M(\log t)^{\frac{\alpha-d+2}{4d}}\left|\int_{\R^d \setminus B_{2M}(t)}
	\nabla v(-y) \bar{\omega}_t({\rm{d}}y)\right|.
\end{split}
\end{equation*}
The integral on the right hand side has zero $\Tilde\P_t$-mean and its variance is
\begin{equation*}
\begin{split}
	\widetilde{\mathbb{V}{\rm ar}}_t \left(\int_{\R^d \setminus B_{2M}(t)}
	\nabla v(-y) \omega({\rm{d}}y)\right)
	&=\nu\int_{\R^d \setminus B_{2M}(t)}
	|\nabla v(-y)|^2 e^{-\rho(\lambda(t))\hat{v}(y)}{\rm{d}}y\\
	&=O\left((\log t)^{\frac{d-2\alpha-2}{d}}\right)
\end{split}
\end{equation*}
due to Lemma~\ref{lem8}-(iii). 
Hence Chebyshev's inequality yields
\begin{equation}
\begin{split}
	&\Tilde\P_t\left(\sup_{x \in B_M(t)}\left|\int_{\R^d \setminus B_{2M}(t)}
	\langle x, \nabla v(-y)\rangle \bar{\omega}_t({\rm{d}}y)\right| 
	> \e (\log t)^{-\frac{\alpha-d+2}{2d}}\right)\\
	&\quad=O\left((\log t)^{-\frac{\alpha+d-2}{2d}}\right)\label{limiter}
\end{split}
\end{equation}
as $t \to \infty$. 
For the second term on the right hand side of~\eqref{B_{2M}(t)^c}, 
we can employ the same argument as above to obtain
\begin{equation*}
\begin{split}
	&\Tilde\P_t\left(\sup_{x \in B_M(t)}\left|\int_{\R^d \setminus B_{2M}(t)}
	\langle x, {\rm Hess}_{v}(-y)x \rangle \bar{\omega}_t({\rm{d}}y)\right|
	> \e (\log t)^{-\frac{\alpha-d+2}{2d}}\right)\\
	&\quad=O\left((\log t)^{-1}\right)
\end{split}
\end{equation*}
Finally, we bound the third term on the right hand side of~\eqref{B_{2M}(t)^c} as 
\begin{equation}
\begin{split}
	&\sup_{x \in B_M(t)}\left|\int_{\R^d \setminus B_{2M}(t)}\int_0^1\frac{(1-\theta)^2}{2}
	\frac{{\rm d}^3}{{\rm d}\theta^3} v(\theta x-y) {\rm d}\theta
	\bar{\omega}_t({\rm{d}}y)\right|\\
	&\quad\le \int_{\R^d \setminus B_{2M}(t)}
	\sup_{x \in B_M(t), \theta \in [0,1]}\left|\frac{{\rm d}^3}{{\rm d}\theta^3} v(\theta x-y) \right|
	\bar{\omega}_t({\rm{d}}y)\\
	&\qquad +2\nu \int_{\R^d \setminus B_{2M}(t)}
	\sup_{x \in B_M(t), \theta \in [0,1]}\left|\frac{{\rm d}^3}{{\rm d}\theta^3} v(\theta x-y)\right|
	e^{-\rho(\lambda(t)) v(y)} {\rm{d}}y.\label{third-term}
\end{split}
\end{equation}
One can easily see that the second term is of $o((\log t)^{-(\alpha-d+2)/2d})$
by using Lemma~\ref{lem8}-(iii) together with~\eqref{3rd-derivative}. 
Furthermore, it also follows that the variance of the first term on the right hand side 
of~\eqref{third-term} is of $O((\log t)^{-(\alpha+d+6)/2d})$. 
Then we can conclude by use of Chebyshev's inequality that 
\begin{equation*}
\begin{split}
	&\Tilde\P_t\left(\int_{\R^d \setminus B_{2M}(t)}
	\sup_{x \in B_M(t), \theta \in [0,1]}\left|\frac{{\rm d}^3}{{\rm d}\theta^3} v(\theta x-y) \right|
	\bar{\omega}_t({\rm{d}}y)
	> \e (\log t)^{-\frac{\alpha-d+2}{2d}}\right)\\
	&\quad =O\left((\log t)^{\frac{\alpha-3d+2}{2d}}\right)=o(1)
\end{split}
\end{equation*}
as $t \to \infty$ and the proof of Lemma~\ref{lem9} is completed.
\qed

\vspace{10pt}
Now we can complete the proof of Proposition~\ref{alpha>2}. Let us define $\omega(z)$ by 
the restriction of $\omega$ on $z+\Lambda_N$ as a measure and introduce events
\begin{align*}
	F_t(z)&=\set{\sup_{x \in z + B_M(t)} 
	|V_{\omega(z)}(x) - Q_t(x-z)| > \frac{\e}{2}(\log t)^{-\frac{\alpha-d+2}{2d}}},\\
	G_t(z)&=\set{\sup_{y \in z+B_M(t)}
        \sum_{\omega_i \not\in z+\Lambda_N}|y-\omega_i|^{-\alpha}
	> \frac{\e}{4} (\log t)^{-\frac{\alpha-d+2}{2d}}}.
\end{align*}
We write $\I_t=\I \cap \{|z|<t(\log t)^{-6}\}$ to simplify notation. 
We are going to show that there exists $z \in \I_t$ for which both $F_t(z)$ and $G_t(z)$ fail to occur. 
To this end, we bound the probability
\begin{equation*}
	\P_{\nu}\left( \bigcap_{z \in \I_t}F_t(z) \cup G_t(z) \right)
	\le \P_{\nu}\left( \bigcap_{z \in \I_t}F_t(z)\right)
	+ \P_{\nu}\left(\bigcup_{z \in \I_t} G_t(z) \right).
\end{equation*}
Note that $\{F_t(z)\}_{z \in \I_t}$ are independent and recall that 
we know 
\begin{equation*}
	\P_{\nu}(G_t(z)) \le \exp\set{-c_7(\log t)^{dN}}
\end{equation*}
from Lemma~\ref{lem6} and that $N > 2$. Moreover, since we have
\begin{equation*}
	F_t(z) \setminus G_t(z) \subset 
	\set{\sup_{x \in z+B_M(t)} |V_{\omega}(x) - Q_t(x-z)| > \frac{\e}{4}
	(\log t)^{-\frac{\alpha-d+2}{2d}}}, 
\end{equation*}
we also know from Lemma~\ref{lem9} that 
\begin{equation*}
\begin{split}
	\P_{\nu}(F_t(z)) &\le \P_{\nu}(F_t(z) \setminus G_t(z)) + \P_{\nu}(G_t(z))\\
	& \le 1-t^{-d}\exp\set{(\log t)^{\delta}}+\exp\set{-c_7(\log t)^{dN}}\\
	& \le 1-\frac{1}{2}t^{-d}\exp\set{(\log t)^{\delta}}. 
\end{split}
\end{equation*}
Combining the above estimates and using $1-x \le e^{-x}$, we obtain 
\begin{equation*}
\begin{split}
	&\P_{\nu}\left( \bigcap_{z \in \I_t}F_t(z) \cup G_t(z) \right)\\
	&\quad \le \left(1-\frac{1}{2}t^{-d}\exp\set{(\log t)^{\delta}}\right)^{t^d(\log t)^{-(6+N)d}}+
	t^d\exp\set{-c_7(\log t)^{dN}}\\
	& \quad \le \exp\set{-\frac{1}{2}e^{(\log t)^{\delta}}(\log t)^{-(6+N)d}}
	+t^d\exp\set{-c_7(\log t)^{dN}}, 
\end{split}	
\end{equation*}
This last expression is summable in $t \in \N$ 
and therefore Borel-Cantelli's lemma tells us that $\P_{\nu}$-almost surely, except for finitely many $t \in \N$, 
there exists $z(t,\omega) \in \I$ with $|z(t,\omega)|<t(\log t)^{-6}$ 
for which both $F_t(z(t,\omega))$ and $G_t(z(t,\omega))$ fail to occur. 

Finally we show that we can take this $z(t,\omega)$ as $x_{2\e,M}(t,\omega)$ for all $t > 0$. 
We can clearly take it as $x_{\e,M}(t,\omega)$ when $t \in \N$. 
Since we have the asymptotic relations 
\begin{align*}
	&(\log (t+1))^{-\frac{\alpha-d}{d}}-(\log t)^{-\frac{\alpha-d}{d}}
	=o\left((\log t)^{-\frac{\alpha-d+2}{2d}}\right),\\
	&(\log (t+1))^{-\frac{\alpha-d+2}{2d}}
	=(\log t)^{-\frac{\alpha-d+2}{2d}}(1+o(1))
\end{align*}
as $t \to \infty$, we can take it as $x_{2\e,M}(t,\omega)$ still in $[t,t+1)$ for large $t$. 
\qed

\subsection{Proof of  Proposition~\ref{local} in the case 
$1< \alpha <2$}\label{remaining}
Note first that $1< \alpha <2$ forces $d=1$ since $\alpha > d$. 
Inspecting the argument in the previous subsection, one can see that we have troubles in \eqref{E_1} 
and \eqref{limiter} in the case $1< \alpha < 2$. 
More precisely, we first need a local limit type theorem to bound the probability in the 
right hand side of \eqref{E_1}. Second, if we have a local limit theorem, the probability in \eqref{E_1}
has order $O((\log t)^{(\alpha-2)/2})$ and then it turns out that \eqref{limiter} is not good enough
when $\alpha \le 3/2$. 

We shall cope with the first problem by proving a local limit theorem around the origin (see Lemma~\ref{lem11}-(i) below). 
For the second problem, we shall bound the fluctuation of 
$\int (V_{\omega}(x)-V_{\omega}(0))\phi_t(x)^2 {\rm d} x$ instead of the potential itself, where $\phi_t(x)$ is 
an approximate eigenfunction (see Lemma~\ref{lem11}-(ii) below). 
Although the latter change prevents us from getting Proposition~\ref{alpha>2}, 
we can show the following lemma which is still sufficient to prove Proposition~\ref{local}. 
\begin{lemma}\label{lem10}
Let $\e>0$ be arbitrarily small. 
There exists $\delta \in (0,1)$ such that for sufficiently large $M>0$ and $t > 0$, 
\begin{equation*}
\begin{split}
	&\P_{\nu}\left( \lambda_{\omega, 1}(B_M(t)) 
	\le \lambda(t) + (q_2+\e)(\log t)^{-\frac{\alpha+1}{2}}\right)\\
	&\quad \ge t^{-1}\exp\set{ (\log t)^{\delta}}.
\end{split}
\end{equation*}
\end{lemma}
This lemma replaces Lemma~\ref{lem9} in the previous subsection and then
Proposition~\ref{local} follows in almost the same way. 
We refrain from repeating the argument and concentrate on proving Lemma~\ref{lem10} in this section.\\

\noindent
\emph{Proof of Lemma~\ref{lem10}} 
By the Rayleigh-Ritz variational formula, the smallest eigenvalue of 
$-\kappa \Delta + V_{\omega}$ in $B_M(t)$ with Dirichlet boundary condition can be 
expressed as 
\begin{equation*}
	\inf_{\stackrel{\phi \in W_0^{1,2}(B_M(t)), }{ \| \phi \|_2=1}}
	\set{\int \kappa \phi'(x)^2 + V_{\omega}(x) \phi(x)^2 {\rm d} x }. 
\end{equation*}
Hence it suffices to find a function $\phi_t \in W_0^{1,2}(B_M(t))$ 
with $\|\phi_t\|_2=1$ satisfying
\begin{equation}
\begin{split}
	&\P_{\nu}\left(\int \kappa \phi'_t(x)^2 
        + V_{\omega}(x)\phi_t(x)^2 {\rm d} x 
	\le \lambda(t) + (q_2+\e)(\log t)^{-\frac{\alpha+1}{2}}\right)\\
	&\quad \ge  t^{-1}\exp\set{(\log t)^{\delta}} \label{to prove}
\end{split}
\end{equation}
for large $t$. 
Let $\bar{\phi}_t$ be the $L^2$-normalized eigenfunction of $-\kappa \Delta + Q_t$ on $\R$ 
corresponding to the smallest eigenvalue $\lambda(t)+q_2(\log t)^{-(\alpha+1)/2}$, that is, 
\begin{equation*}
	\bar{\phi}_t(x)=\left(\frac{q_2}{\kappa \pi}\right)^{\frac{1}{4}}(\log t)^{-\frac{\alpha+1}{8}}
	\exp\set{-\frac{q_2}{2\kappa}(\log t)^{-\frac{\alpha+1}{2}}|x|^2}.
\end{equation*}
We define $\phi_t \in W_0^{1,2}(B_M(t))$ by 
\begin{equation*}
	\phi_t= c_M(t)(\bar{\phi}_t-\bar{\phi}_t(M(\log t)^{\frac{\alpha+1}{4}})) 1_{B_M(t)},
\end{equation*}
where $c_M(t)$ is chosen so that $\|\phi_t\|_2=1$. 
One can easily check that $c_M(t) \to 1$ as $M \to \infty$. 
Note also that we have
\begin{equation*}
\begin{split}
	&\int \kappa \phi_t'(x)^2 + V_{\omega}(x) \phi_t(x)^2 {\rm d} x \\
	&\quad \le c_M(t)^{-2} \int_{B_M(t)} 
	\kappa \bar{\phi}'_t(x)^2 + V_{\omega}(x) \bar{\phi}_t(x)^2 {\rm d} x.
\end{split}
\end{equation*}
Hence it suffices to show \eqref{to prove} with $\phi_t$ replaced  by $\bar{\phi}_t\cdot 1_{B_M(t)}$ 
when $M$ is sufficiently large. 
Let us introduce the events
\begin{align*}
	&E_1=\set{(\log t)^{\frac{\alpha+1}{2}}(V_{\omega}(0) - \lambda(t)) \in (\e/4, \e/2)},\\
	&E_2=\set{\int_{B_M(t)} ([V_{\omega}(x)-V_{\omega}(0)]-[Q_t(x)-\lambda(t)])\bar{\phi}_t(x)^2 {\rm d} x 
	\ge \frac{\e}{2} (\log t)^{-\frac{\alpha+1}{2}}}.
\end{align*}
Since we have 
\begin{equation*}
	\int_{B_M(t)} \kappa \bar{\phi}'_t(x)^2 + 
        Q_t(x) \bar{\phi}_t(x)^2 {\rm d} x 
	=\lambda(t)+(q_2+o(1))(\log t)^{-\frac{\alpha+1}{2}}
\end{equation*}
as $M \to \infty$ by the definition of $\bar{\phi}_t$, we see
\begin{equation*}
\begin{split}
	\P_{\nu}&\left(\int_{B_M(t)} 
	\kappa \bar{\phi}'_t(x)^2 + V_{\omega}(x)\bar{\phi}_t(x)^2 {\rm d} x
	\le \lambda(t) + (q_2+2\e)(\log t)^{-\frac{\alpha+1}{2}}\right)\\
	&= \P_{\nu} \left( \int_{B_M(t)} (V_{\omega}(x)-Q_t(x))\bar{\phi}_t(x)^2 {\rm d} x 
	\le \e (\log t)^{-\frac{\alpha+1}{2}} \right)\\
	&\ge \P_{\nu}(E_1 \setminus E_2)
\end{split}
\end{equation*}
for sufficiently large $M$. 
Recalling the definition of transformed measure $\Tilde{\P}_t$, \eqref{rho(lambda)}, and \eqref{dlogt}, 
we have
\begin{equation*}
	\P_{\nu}(E_1 \setminus E_2)=\exp\set{-\log t+o(1)+\frac{\e}{4}\rho(e)(\log t)^{\frac{\alpha-1}{2}}}
	\Tilde{\P}_t(E_1 \setminus E_2)
\end{equation*}
as $t \to \infty$. 
\begin{lemma}\label{lem11}
\begin{enumerate}
\item{When $d=1$, there exists a constant $c_9(\nu, \alpha)>0$ such that 
\begin{equation}
	\lim_{t \to \infty}
        \Tilde{\P}_t \Bigl((\log t)^{\frac{2\alpha-1}{2}}(V_{\omega}(0) - \lambda(t)) \in 
	(0, a)\Bigr) = c_9 a(1+o(1)) \label{lem11-1}
\end{equation}
as $a \downarrow 0$.}
\item{Let $d=1$ and $\bar{\phi}_t$ be as above. Then
\begin{equation}
	\widetilde{\mathbb{V}{\rm ar}}_t\left(\int_{B_M(t)} 
	(V_{\omega}(x)-V_{\omega}(0))\bar{\phi}_t(x)^2 {\rm d} x \right)
	=O((\log t)^{-\alpha-2}) \label{lem11-2}
\end{equation}
as $t \to \infty$. }
\end{enumerate}
\end{lemma}
We defer the proof of this lemma and finish the proof of Lemma~\ref{lem10} first. 
By using \eqref{lem11-1}, we obtain 
\begin{equation}
	\Tilde{\P}_t (E_1) = 
        c_9\frac{\e}{4} (\log t)^{\frac{\alpha-2}{2}}(1+o(1))\label{lem11-1'} 
\end{equation}
as $t \to \infty$. To bound $\Tilde{\P}_t (E_2)$, we first
replace $Q_t(x)-\lambda(t)$ in $E_2$ by 
$\Tilde{\E}_t[V_{\omega}(x)-V_{\omega}(0)]$ using Lemma~\ref{lem7}-(ii). 
Then, by \eqref{lem11-2} and Chebyshev's inequality, it follows that
\begin{equation}
	\Tilde{\P}_t (E_2)=O\left((\log t)^{-1}\right)
	= o\left((\log t)^{\frac{\alpha-2}{2}}\right) 
        \quad {\rm as} \quad  t \to \infty\label{lem11-2'}
\end{equation}
Combining \eqref{lem11-1'} and \eqref{lem11-2'}, we obtain 
\begin{equation*}
	\Tilde{\P}_t(E_1 \setminus E_2) 
        \ge \frac{c_9\e}{8}(\log t)^{\frac{\alpha-2}2}
\end{equation*}
and the proof of Lemma \ref{lem10} is completed. 
\qed

\vspace{10pt}\noindent
\emph{Proof of Lemma \ref{lem11} {\upshape (i)}}
Let 
\begin{align*}
	f_{\rho}(\theta) 
	& = \exp\set{\nu \int \Bigr(e^{i \theta\rho^{\frac{2\alpha-1}{2\alpha}}\hat{v}(y)}
	-1 - i \theta \rho^{\frac{2\alpha-1}{2\alpha}}\hat{v}(y)\Bigr) 
	e^{-\rho\hat{v}(y)} {\rm d} y}, \\
	g_{\rho}(\theta)
	&=\exp\set{-\frac12 \nu \theta^2  \rho^{\frac{2\alpha-1}{\alpha}}
	 \int \hat{v}(y)^2 e^{-\rho\hat{v}(y)} {\rm d} y}.
\end{align*}
Note that $f_{\rho(\lambda(t))}$ is the characteristic function of the law $\mu_t$ of 
$\rho(\lambda(t))^{(2\alpha-1)/2\alpha}(V_{\omega}(0)-\lambda(t))$ under $\Tilde{\P}_t$
by \eqref{characteristic} 
and that $g_{\rho(\lambda(t))}$ is that of a Gaussian measure $\nu_t$ whose variance 
converges to $\nu \int |\eta|^{-2\alpha} e^{-|\eta|^{-\alpha}} {\rm d} \eta$ as $t \to \infty$. 
Therefore it suffices to prove that 
$a^{-1}|\mu_t([0,a])-\nu_t([0,a])|\to 0$ as 
$t \to \infty$ uniformly in $a>0$.

By using L\'{e}vy's inversion formula and the fact that $(e^{i x}-1)/x$ is bounded for $x \in \R$, 
we have 
\begin{equation}
\begin{split}
	&\frac{1}{a} | \mu_t([0,a])-\nu_t([0,a]) | \\
	&\quad = \lim_{T \to \infty} \left| \int_{-T}^T \frac{1-e^{-i a\theta}}{i a\theta}
	(f_{\rho}(\theta)-g_{\rho}(\theta)) {\rm d} \theta \right| \\
	&\quad \le \int_{|\theta| \le \e \rho^{1/3\alpha}} 
	|f_{\rho}(\theta)-g_{\rho}(\theta)| {\rm d} \theta 
	+\int_{ |\theta| > \e \rho^{1/3\alpha}}
	|g_{\rho}(\theta)| {\rm d} \theta \\
	&\hspace{30pt} 
	+\int_{\e \rho^{1/3\alpha}< |\theta| \le \e \rho^{1/2\alpha}}
	|f_{\rho}(\theta)| {\rm d} \theta 
	+\int_{|\theta| > \e \rho^{1/2\alpha}} 
	|f_{\rho}(\theta)| {\rm d} \theta \\
	& \quad =: I_1 + I_2 + I_3 + I_4,\label{inversion}
\end{split}
\end{equation}
where $\e >0$ is a small constant which will be chosen later. 
Note first that $I_2 \to 0 $ as $\rho \to \infty$. 
To bound $I_1$, we rewrite the integrand as
\begin{equation*}
\begin{split}
	&|f_{\rho}(\theta)-g_{\rho}(\theta)| \\
	&\quad = g_{\rho}(\theta) 
	\biggl|\exp\biggl\{\nu \int \Bigr(e^{i \theta\rho^{\frac{2\alpha-1}{2\alpha}}\hat{v}(y)}-1
	- i \theta \rho^{\frac{2\alpha-1}{2\alpha}}\hat{v}(y) \\
	&\hspace{130pt} + \frac12\theta^2  \rho^{\frac{2\alpha-1}{\alpha}}\hat{v}(y)^2 \Bigr) 
	e^{-\rho\hat{v}(y)} {\rm d} y\biggr\} -1 \biggr|. 
\end{split}
\end{equation*}
By using the bound 
\begin{equation*}
	\left| e^{i x}-1- i x+ \frac{1}{2}x^2 \right| 
	= \left| \frac{1}{3!}\int_0^x (x-s)^2 e^{i s} {\rm d} s \right| 
	\le \frac{1}{6} |x|^3, 
\end{equation*}
and change of variable, we obtain 
\begin{equation*}
\begin{split}
	&\int \Bigr(e^{i \theta\rho^{\frac{2\alpha-1}{2\alpha}}\hat{v}(y)}-1
	- i \theta \rho^{\frac{2\alpha-1}{2\alpha}}\hat{v}(y) 
	+ \frac12\theta^2  \rho^{\frac{2\alpha-1}{\alpha}}\hat{v}(y)^2 \Bigr) 
	e^{-\rho\hat{v}(y)} {\rm d} y \\
	&\quad \le |\theta|^3  \rho^{\frac{6\alpha-3}{2\alpha}} 
	\int \hat{v}(y)^3 e^{-\rho\hat{v}(y)} {\rm d} y \\
	&\quad = O\left( |\theta|^3  \rho^{-\frac{1}{\alpha}} \right). 
\end{split}
\end{equation*}
Hence if $|\theta| \le \e \rho^{1/3\alpha}$ and $\e$ is sufficiently small, 
then an elementary inequality $|e^z-1| < 2z$ which holds for small $z \in \C$ yields 
\begin{equation*}
	|f_{\rho}(\theta)-g_{\rho}(\theta)| 
	= O\left(g_{\rho}(\theta) |\theta|^3  \rho^{-\frac{1}{\alpha}}\right).
\end{equation*}
This shows that 
\begin{equation*}
	I_1 \le O\left( \int g_{\rho}(\theta) |\theta|^3  \rho^{-\frac{1}{\alpha}} {\rm d} \theta\right) 
	=o(1) 
\end{equation*}
as $\rho \to \infty$. For larger $|\theta|$, we use 
\begin{equation*}
	|f_{\rho}(\theta)|=
	\exp\set{-\nu \int (1-\cos(\theta\rho^{\frac{2\alpha-1}{2\alpha}}\hat{v}(y)))
	e^{-\rho \hat{v}(y)} {\rm d} y}.
\end{equation*}
If $\e \rho^{1/3\alpha} < |\theta| < \e \rho^{1/2\alpha}$ and $\e$ is sufficiently small, 
then by using $1-\cos x \ge x^2/4$ for small $x \in \R$, we have
\begin{equation*}
\begin{split}
	&\int (1-\cos(\theta\rho^{\frac{2\alpha-1}{2\alpha}}\hat{v}(y))) e^{-\rho \hat{v}(y)} {\rm d} y\\
	&\quad \ge \int_{ \rho^{1/\alpha}}^{2 \rho^{1/\alpha}}
	\frac{\theta^2}{4} \rho^{\frac{2\alpha-1}{\alpha}}|y|^{-2\alpha} e^{-\rho \hat{v}(y)}{\rm d} y \\ 
	&\quad \ge 2^{-2\alpha-2}e^{-1}\theta^2.
\end{split}
\end{equation*}
It follows from this that 
\begin{equation*}
	I_2 \le \int_{ |\theta| > \e \rho^{1/3\alpha}}
	\exp\set{-2^{-2\alpha-2}e^{-1} \theta^2} {\rm d} \theta \xrightarrow{\rho \to \infty} 0.
\end{equation*}
Finally if $|\theta| \ge \e \rho^{1/2\alpha}$, then 
\begin{equation*}
\begin{split}
	&\int (1-\cos(\theta\rho^{\frac{2\alpha-1}{2\alpha}}\hat{v}(y))) e^{-\rho \hat{v}(y)} {\rm d} y\\
	& \quad \ge \int_{|\theta|^{1/\alpha} \rho^{(2\alpha-1)/2\alpha^2}}
	^{2 |\theta|^{1/\alpha} \rho^{(2\alpha-1)/2\alpha^2}} 
	(1-\cos(\theta\rho^{\frac{2\alpha-1}{2\alpha}}|y|^{-\alpha})) e^{-\rho |y|^{-\alpha}}{\rm d} y \\ 
	& \quad \ge (1-\cos 2^{-\alpha})e^{-\e^{-\alpha}}
	\rho^{\frac{2\alpha-1}{2\alpha^2}}\theta^{\frac{1}{\alpha}},
\end{split}	
\end{equation*}
and hence 
\begin{equation*}
	I_4 \le \int_{|\theta| > \e \rho^{1/2\alpha}} 
	\exp\set{- (1-\cos 2^{-\alpha})e^{-\e^{-\alpha}} \rho^{\frac{2\alpha-1}{2\alpha^2}}
	\theta^{\frac{1}{\alpha}}} {\rm d} \theta 
	\xrightarrow{\rho \to \infty} 0.
\end{equation*}
Coming back to \eqref{inversion}, we see that 
\begin{equation*}
	\frac{1}{a} | \mu_t([0,a])-\nu_t([0,a]) | \rightarrow 0 
\end{equation*}
as $t \to \infty$ uniformly in $a>0$ and the claim follows. 
\qed

\vspace{10pt}\noindent
\emph{Proof of Lemma \ref{lem11} {\upshape (ii)}}
We use \eqref{variance} to obtain
\begin{equation*}
\begin{split}
	&\widetilde{\mathbb{V}{\rm ar}}_t \left(\int_{B_M(t)} 
	(V_{\omega}(x)-V_{\omega}(0))\bar{\phi}_t(x)^2 {\rm d} x \right)\\
	&\quad=\nu\int \set{\int_{B_M(t)} 
	(\hat{v}(x-y)-\hat{v}(-y))\bar{\phi}_t(x)^2 {\rm d} x}^2 e^{-\rho(\lambda(t))\hat{v}(-y)} {\rm d} y.
\end{split}
\end{equation*}
By using the facts that $\hat{v}$ is bounded and $\|\bar{\phi}_t\|_2=1$ together with Lemma~\ref{lem8}-(i),  
we see that the integral over the region $y \in B_{2M}(t)$ makes only negligible
contribution in the right hand side.
On the region $y \not \in B_{2M}(t)$, we may replace $\hat{v}$ by $v$ 
and it follows by Taylor's theorem that
\begin{equation*}
\begin{split}
	&\int_{B_M(t)} (v(x-y)-v(-y))\bar{\phi}_t(x)^2 {\rm d} x\\
	&\quad = \int_{B_M(t)} \left(v'(-y) x  + 
	\int_0^1(1-\theta)\frac{{\rm d}^2}{{\rm d} \theta ^2}v(\theta x-y)  {\rm d} \theta \right)
	\bar{\phi}_t(x)^2 {\rm d} x\\
	&\quad = v'(-y) \int_{B_M(t)} x \bar{\phi}_t(x)^2 {\rm d} x \\ 
	&\hspace{40pt} + \int_{B_M(t)} \int_0^1 
	(1-\theta)\frac{{\rm d}^2}{{\rm d} \theta ^2}v(\theta x-y) {\rm d} \theta \,\bar{\phi}_t(x)^2 {\rm d} x. 
\end{split}
\end{equation*}
The first term in the last line vanishes since $x \bar{\phi}_t(x)^2$ is 
symmetric about the origin. 
The second term is of order
$$
O\left(|y|^{-\alpha-2}\int_{B_M(t)}|x|^2\bar{\phi}_t(x)^2{\rm d}x\right)
=O\left(|y|^{-\alpha-2}(\log t)^{\frac{\alpha+1}{2}}\right)
$$ 
uniformly in $y \not\in B_{2M}(t)$.
Therefore, we arrive at 
\begin{equation*}
\begin{split}
 	&\widetilde{\mathbb{V}{\rm ar}}_t \left(\int_{B_M(t)}
	 (V_{\omega}(x)-V_{\omega}(0))\bar{\phi}_t(x)^2 {\rm d} x \right)\\
	&\quad=O\left( (\log t)^{\alpha+1}
         \int_{\R^d\setminus B_{2M}(t)} 
       |y|^{-2\alpha-4} e^{-\rho(\lambda(t))|y|^{-\alpha}} 
       {\rm d} y \right)
\end{split}
\end{equation*}
and the result follows by applying Lemma~\ref{lem8}-(iii). 
\qed

\section{Lower bound on the integrated density of states}\label{Lifshiz-lower}
In this section, we prove the lower bound of Theorem~\ref{Lifshitz}. 
Since there seems to be no exact bound like \eqref{exact}, 
we cannot derive the second order asymptotics from Theorem~\ref{annealed}. 
We instead use \eqref{IDS-lower} to reduce the problem to the estimate of the 
principal eigenvalue in a finite box and then use Proposition~\ref{local}. \\

\noindent
\emph{Proof of the lower bound in Theorem~\ref{Lifshitz}}
Let us fix $\e > 0$ arbitrarily small and take 
$R=t$ and $\lambda=\lambda(t,-\e)$ (see \eqref{lambda} for the definition) 
in \eqref{IDS-lower} to obtain 
\begin{equation*}
	N(\lambda(t,-\e)) \ge (2t)^{-d}\P_{\nu} \left( \lambda_{\omega,\, 1}\bigl((-t, t)^d\bigr) \le \lambda(t,-\e) \right). 
\end{equation*}
It is straightforward to check that 
\begin{equation*}
	t^{-d} \ge \exp\set{ -l_1 \lambda(t,-\e)^{-\frac{d}{\alpha-d}} 
	- l_2 \lambda(t,-\e)^{-\frac{\alpha+d-2}{2(\alpha-d)}}}
\end{equation*}
for any $\e>0$ when $t$ is sufficiently large. 
Therefore, it suffices to prove that
\begin{equation*}
	\P_{\nu} \left( \lambda_{\omega,\, 1}\bigl((-t, t)^d\bigr) \le \lambda(t,-\e) \right) 
	\to 1
\end{equation*}
as $t \to \infty$. This convergence follows from Proposition~\ref{local}. 
Indeed, it implies that the probability of having a ball 
\begin{equation*}
	B \left(x_{\e,M}(t,\omega), M(\log t)^{\frac{\alpha-d+2}{4d}}\right), \;
	| x_{\e,M}(t,\omega) | \le t(\log t)^{-6}
\end{equation*}
such that 
\begin{equation*}
	\lambda_{\omega,\, 1}\left(B_{\e,M}(t,\omega)\right)
	\le q_1 (\log t)^{-\frac{\alpha-d}{d}}+
	(q_2+\e)(\log t)^{-\frac{\alpha-d+2}{2d}}
\end{equation*}
approaches to 1 as $t \to \infty$ and, needless to say, 
$ \lambda_{\omega,\, 1}((-t, t)^d) \le \lambda_{\omega,\, 1}(B_{\e,M}(t,\omega))$. 
\qed

\section{Appendix}\label{appendix}
We collect some formulae for Poisson point process 
which we use in this paper. 
\begin{proposition}
 Let $(\omega, \P_{m})$ be the Poisson point process with intensity 
 $m({\rm d} x)$ being a positive Radon measure. 
 \begin{enumerate}
 \item If $f$ is a sign definite Borel function, 
 \begin{equation}
	\E_m\left[ \exp\set{\int f(y) \,\omega({\rm d} y)} \right]
	=\exp\set{\int (e^{f(y)}-1)\,m({\rm d} y)}.\label{Laplace}
 \end{equation}
 \item If $f$ is an $m$-integrable real valued function, then
 \begin{equation}
	\E_m\left[ \exp\set{i \int f(y) \,\omega({\rm d} y)} \right]
	=\exp\set{\int (e^{i f(y)}-1)\,m({\rm d} y)}.\label{characteristic}
 \end{equation}
 \item If $f$ is an $m$-integrable function, then 
\begin{equation}
	\E_m \left[\int f(y) \,\omega ({\rm d} y) \right]
	=\int f(y) \,m( {\rm d} y)\label{expectation}
\end{equation}
 \item If both $f$ and $f^2$ are $m$-integrable, then
\begin{equation}
	{\mathbb{V}{\rm ar}}_m \left(\int f(y) \,\omega ({\rm d} y) \right)
	=\int f(y)^2 \,m( {\rm d} y).\label{variance}
\end{equation}
 \end{enumerate}
\end{proposition}
\emph{Proof}
The first two formulae are Lemma~10.2 (p.178) in 
Kallenberg~\cite{Kal02}.
In fact, the formula~\eqref{Laplace} is proved only for 
non-positive functions there but the argument can easily 
be adapted to non-negative case. 

The latter two assertions follow by differentiating
~\eqref{characteristic} if $f \in C_c(\R^d)$.
Then they can be generalized as stated above by approximation.
\qed

\section*{Acknowledgements}
This work emerged from a discussion with professor Wolfgang K\"{o}nig. 
The author thanks him for useful discussions about this problem. 
Part of this work was done when the author was visiting Max-Planck-Institut 
f\"ur Mathematik in den Naturwissenschaften in Leipzig and he would like to 
thank the warm hospitality. 

\newcommand{\noop}[1]{}

\end{document}